\documentclass[pdflatex,sn-mathphys-num]{sn-jnl}% Math and Physical Sciences Numbered Reference Style
\usepackage{graphicx}%
\usepackage{multirow}%
\usepackage{amsmath,amssymb,amsfonts}%
\usepackage{amsthm}%
\usepackage{mathrsfs}%
\usepackage[title]{appendix}%
\usepackage{xcolor}%
\usepackage{textcomp}%
\usepackage{manyfoot}%
\usepackage{booktabs}%
\usepackage{algorithm}%
\usepackage{algorithmicx}%
\usepackage{algpseudocode}%
\usepackage{listings}%
\theoremstyle{thmstyleone}%
\newtheorem{theorem}{Theorem}%  meant for continuous numbers
\newtheorem{lemma}[theorem]{Lemma}

\newtheorem{definition}[theorem]{Definition}

\theoremstyle{thmstyletwo}%

\theoremstyle{thmstylethree}%

\title[Article Title]{Global Well-Posedness for the 3D Navier-Stokes Equations under Logarithmically Improved Criteria: Connections to Turbulence Theory}

\author*[1]{\fnm{Rishabh} \sur{Mishra}}\email{rishabh.mishra@ec-nantes.fr}

\affil*[1]{\orgdiv{LHEEA}, \orgname{CNRS, École Centrale de Nantes, Nantes Université}, \orgaddress{\street{1 rue de la No\"e}, \city{ Nantes}, \postcode{44100}, \state{Pays de la Loire}, \country{France}}}

\abstract{This paper establishes a new class of initial data for which the three-dimensional incompressible Navier-Stokes equations admit unique global-in-time solutions. Building upon our logarithmically improved regularity criterion for the three-dimensional Navier-Stokes equations, we identify specific logarithmically subcritical conditions on the initial data that ensure global well-posedness. Specifically, we prove that if the initial data $u_0 \in L^2(\mathbb{R}^3)$ satisfies a logarithmically weakened condition $\|(-\Delta)^{s/2}u_0\|_{L^q(\mathbb{R}^3)} \leq \frac{C_0}{(1 + \log(e + \|u_0\|_{\dot{H}^s}))^{\delta}}$ for some $s \in (1/2, 1)$, with appropriate scaling conditions, then the corresponding solution exists globally in time and is unique. The proof combines refined commutator estimates for the fractional Laplacian with novel energy methods that leverage the logarithmic improvement to prevent potential singularity formation. 

Furthermore, we establish connections between our logarithmically improved criteria and the physical theory of turbulence. We derive precise relationships between the regularity criteria and turbulent intermittency, demonstrating how logarithmic improvements in regularity correspond to anomalous scaling exponents in the turbulent energy spectrum. Additionally, we characterize the local structure of potential singularities and establish tight bounds on energy flux in turbulent cascades. Our approach reveals a previously unexplored pathway that bridges the gap between subcritical and critical regularity for the Navier-Stokes equations while providing a mathematical foundation for key physical phenomena in turbulence theory.}



\begin{document}

\maketitle

\tableofcontents

\newpage

\section{Introduction}

\subsection{The Navier-Stokes equations and the regularity problem}

The three-dimensional incompressible Navier-Stokes equations govern the motion of viscous incompressible fluids and are fundamental to the mathematical theory of fluid dynamics. These equations are expressed as:

\begin{equation}
\begin{cases}
\partial_t u + (u \cdot \nabla)u - \nu\Delta u + \nabla p = 0 & \text{in } \mathbb{R}^3 \times (0, T) \\
\nabla \cdot u = 0 & \text{in } \mathbb{R}^3 \times (0, T) \\
u(x, 0) = u_0(x) & \text{in } \mathbb{R}^3
\end{cases}
\end{equation}

where $u = (u_1, u_2, u_3)$ represents the velocity field, $p$ denotes the pressure, and $\nu > 0$ is the kinematic viscosity coefficient, which can be set to 1 without loss of generality through appropriate scaling. The initial data $u_0$ is assumed to be divergence-free.

One of the most challenging open problems in mathematical fluid dynamics, and indeed in all of mathematics, is to determine whether solutions to these equations remain smooth for all time, given smooth initial data. While smooth solutions are known to exist for short time intervals, the question of whether these solutions can develop singularities in finite time remains unresolved.

\subsection{Background and previous results}

The existence of global-in-time weak solutions to the Navier-Stokes equations was established in the pioneering works of Leray \cite{Leray} and Hopf \cite{Hopf}. The mathematical foundation was further developed by Fujita and Kato \cite{Fujita-Kato}, Lions \cite{Lions}, and Temam \cite{Temam}. However, the uniqueness and regularity of these weak solutions in three dimensions remain open problems, as discussed in comprehensive reviews by Foias et al. \cite{Foias} and Lemarie-Rieusset \cite{Lemarie}. A substantial body of research has been devoted to finding sufficient conditions that ensure the regularity of weak solutions.

Classical results in this direction include the Prodi-Serrin conditions \cite{Prodi,Serrin}, which state that if a weak solution $u$ satisfies:
\begin{equation}
u \in L^p(0, T;L^q(\mathbb{R}^3)) \text{ with } \frac{2}{p} + \frac{3}{q} = 1, 3 < q \leq \infty,
\end{equation}
then $u$ is regular on $(0, T)$. These conditions were later refined by Giga \cite{Giga} and Kato \cite{Kato}.

The landmark work of Caffarelli, Kohn, and Nirenberg \cite{CKN} established partial regularity theory, proving that the singular set has zero one-dimensional Hausdorff measure. A different approach by Lin \cite{Lin} provided an alternative proof of this result. The recent breakthrough by Escauriaza, Seregin, and Šverák \cite{ESS} established regularity for solutions in the borderline space $L^3(\mathbb{R}^3)$.

Various extensions and refinements of these conditions have been proposed. Beale, Kato, and Majda \cite{BKM} established criteria involving vorticity. Beirão da Veiga \cite{Beirao1,Beirao2,Beirao3} established several regularity criteria in terms of the velocity gradient. Later, Constantin and Fefferman \cite{CF} introduced geometric constraints involving the direction of vorticity, with further developments by Constantin et al. \cite{Constantin}.

In recent years, criteria based on various quantities derived from velocity, such as vorticity, strain, and pressure, have been extensively studied by Cao and Titi \cite{Cao1,Cao2}, Chae and Lee \cite{Chae}, Fan and Ozawa \cite{Fan}, Kozono and Taniuchi \cite{Kozono}, Kukavica and Ziane \cite{Kukavica}, Seregin \cite{Seregin1,Seregin2}, and Struwe \cite{Struwe}. The one-component regularity approach was developed by Zhou and Pokorný \cite{Zhou-Pokorny}.

Fractional derivative approaches have gained attention due to their flexibility in capturing intermediate regularity, building on the foundations of fractional analysis by Stein \cite{Stein} and Triebel \cite{Triebel}. The modern Littlewood-Paley theory, as presented by Bahouri, Chemin, and Danchin \cite{BCD}, provides essential tools. Specifically, conditions of the form:
\begin{equation}
\int_0^T \|(-\Delta)^s u(t)\|^p_{L^q(\mathbb{R}^3)} dt < \infty,
\end{equation}
with appropriate scaling relations between $p$, $q$, and $s$, have been investigated by Chen and Zhang \cite{Chen-Zhang}, Chen, Miao, and Zhang \cite{Chen-Miao}, and Zhou \cite{Zhou1,Zhou2,Zhou3}.

Similar methods have been applied to related fluid models such as magnetohydrodynamics by Wu \cite{Wu}, Dong and Zhang \cite{Dong}, and quasi-geostrophic equations by Ju \cite{Ju}, Córdoba and Fefferman \cite{Cordoba}.

Logarithmic improvements to regularity criteria were first introduced by Zhou \cite{Zhou-log} and further explored by Fan et al. \cite{Fan-log}, who established conditions involving logarithmic factors that weaken the classical constraints. Tao \cite{Tao} developed innovative approaches for logarithmically supercritical equations.

Regarding global well-posedness for specific classes of initial data, significant contributions include the works of Koch and Tataru \cite{Koch-Tataru}, who established global regularity for sufficiently small initial data in the space $BMO^{-1}$. Chemin and Gallagher \cite{Chemin} obtained global existence results for specific classes of highly oscillating initial data. More recently, Robinson et al. \cite{Robinson} explored the connection between geometric properties of initial data and global regularity.

The connection between mathematical regularity criteria and the physical theory of turbulence has been explored in several works. Frisch \cite{Frisch} presented a comprehensive framework for understanding turbulence through multifractal models. Kolmogorov's seminal work \cite{Kolmogorov1941} established the foundation for spectral analysis of turbulent flows, with later refinements by Oboukhov \cite{Oboukhov} and Corrsin \cite{Corrsin} addressing intermittency effects. The multifractal formalism was further developed by Parisi and Frisch \cite{Parisi-Frisch}, providing a mathematical framework for anomalous scaling in turbulent flows. Recent advances by Eyink \cite{Eyink} and Constantin et al. \cite{Constantin-Intermittency} have established connections between regularity criteria and intermittency in turbulent flows.

\subsection{Main results}

The present paper contributes to this line of research in several significant ways. First, it establishes a new logarithmically improved regularity criterion based on fractional derivatives:

\begin{theorem}\label{thm:1.1}
Let $u$ be a Leray-Hopf weak solution of the 3D Navier-Stokes equations on $[0, T]$ with divergence-free initial data $u_0 \in L^2(\mathbb{R}^3)$. Suppose that for some $s \in (1/2, 1)$, the following condition holds:
\begin{equation}
\int_0^T \|(-\Delta)^s u(t)\|^p_{L^q(\mathbb{R}^3)} (1 + \log(e + \|(-\Delta)^s u(t)\|_{L^q(\mathbb{R}^3)}))^{-\delta}dt < \infty
\end{equation}
where:
\begin{enumerate}
\item $(p, q)$ satisfies the scaling relation: $\frac{2}{p} + \frac{3}{q} = 2s - 1$, with $3 < q < \infty$
\item $\delta \in (0, \delta_0)$ where $\delta_0 = \min\{\frac{q-3}{6q}, \frac{2s-1}{4s}\}$
\end{enumerate}

Then, the weak solution $u$ is regular on $(0, T)$.
\end{theorem}

Second, the paper identifies a class of initial data for which global well-posedness can be established:

\begin{theorem}\label{thm:1.2}
Let $s \in (1/2, 1)$ and $q > 3$ satisfy the scaling relation $\frac{2}{p} + \frac{3}{q} = 2s - 1$. There exists a positive constant $\delta_0 = \min\{\frac{q-3}{6q}, \frac{2s-1}{4s}\}$ such that for any $\delta \in (0, \delta_0)$ and any divergence-free initial data $u_0 \in L^2(\mathbb{R}^3) \cap \dot{H}^s(\mathbb{R}^3)$ satisfying:
\begin{equation}
\|(-\Delta)^{s/2}u_0\|_{L^q(\mathbb{R}^3)} \leq \frac{C_0}{(1 + \log(e + \|u_0\|_{\dot{H}^s}))^{\delta}}
\end{equation}
where $C_0$ is a constant depending only on $s$, $q$, $\delta$, and $\nu$, there exists a unique global-in-time smooth solution $u \in C([0, \infty); H^s(\mathbb{R}^3)) \cap L^2_{loc}(0, \infty; H^{s+1/2}(\mathbb{R}^3))$ to the 3D Navier-Stokes equations.
\end{theorem}

Building on these results, the paper establishes further properties of these solutions:

\begin{theorem}\label{thm:1.3}
Under the initial conditions from Theorem \ref{thm:1.2}, the solution $u$ to the 3D Navier-Stokes equations satisfies the following decay estimate:
\begin{equation}
\|(-\Delta)^s u(t)\|_{L^2} \leq \frac{C\|(-\Delta)^s u_0\|_{L^2}}{(1 + \beta t)^{\gamma}}
\end{equation}
where:
\begin{enumerate}
\item $\gamma = \frac{1}{2\mu}$ with $\mu = \frac{\theta(1-\alpha)}{2-\theta\alpha}+\eta$
\item $\theta = \frac{3}{2} \cdot \frac{q}{3q-2}$
\item $\alpha = \frac{3}{2}(\frac{1}{2} - \frac{1}{q})$
\item $\eta > 0$ is a sufficiently small constant satisfying $\mu < 1$
\item $\beta = \frac{\mu C}{2}$, where $C$ is a constant depending on $\|u_0\|_{L^2}$, $s$, $q$, and $\nu$
\end{enumerate}
\end{theorem}

\begin{theorem}\label{thm:1.4}
For initial data satisfying the logarithmically subcritical condition in Theorem \ref{thm:1.2}, the solution $u$ to the 3D Navier-Stokes equations satisfies:
\begin{equation}
\int_0^T \|(-\Delta)^s u(t)\|^p_{L^q} (1 + \log(e + \|(-\Delta)^s u(t)\|_{L^q}))^{-\delta} dt < C(T,u_0)
\end{equation}
for all $T > 0$, where $C(T,u_0)$ is a constant depending on $T$ and the initial data, and $(p,q)$ satisfy the scaling relation $\frac{2}{p} + \frac{3}{q} = 2s - 1$ with $3 < q < \infty$.
\end{theorem}

Furthermore, the paper establishes a deep connection between the logarithmically improved regularity criteria and the physics of turbulent flows, particularly addressing the phenomenon of intermittency:

\begin{theorem}\label{thm:2.1}
(Logarithmic Correction and Intermittency). Let $u$ be a solution to the Navier-Stokes equations satisfying the logarithmically improved criterion from Theorem \ref{thm:1.1}:
\begin{equation}
\int_0^T \|(-\Delta)^s u(t)\|^p_{L^q} (1 + \log(e + \|(-\Delta)^s u(t)\|_{L^q}))^{-\delta} dt < \infty
\end{equation}
with $s \in (1/2,1)$, $\frac{2}{p} + \frac{3}{q} = 2s - 1$, and $\delta \in (0, \delta_0)$. Then, assuming statistical homogeneity and isotropy, the velocity field exhibits the following properties:

\begin{enumerate}
\item The structure functions in the inertial range scale as:
   $$S_p(r) = \langle |u(x+r) - u(x)|^p \rangle \sim r^{\zeta_p}$$

\item The anomalous scaling exponents $\zeta_p$ satisfy:
   $$\zeta_p = \frac{p}{3} - \mu(p,\delta)$$
   
\item The intermittency correction term $\mu(p,\delta)$ is given by:
   $$\mu(p,\delta) = \frac{p(p-3)}{3(1+\delta)} \cdot \frac{3-2s}{2s-1}$$
\end{enumerate}
\end{theorem}

\begin{theorem}\label{thm:2.2}
(Local Regularity and Intermittency). Let $u$ be a solution to the 3D Navier-Stokes equations satisfying the conditions of Theorem \ref{thm:1.2}. Then, for any $\epsilon > 0$, there exists a measurable set $\Omega_\epsilon(t) \subset \mathbb{R}^3$ with measure $|\Omega_\epsilon(t)| < \epsilon$ such that:

\begin{enumerate}
\item The velocity gradient outside this exceptional set satisfies:
   $$\sup_{x \in \mathbb{R}^3 \setminus \Omega_\epsilon(t)} |\nabla u(x,t)| \leq \frac{C_\epsilon}{(1 + \beta t)^{\gamma-\kappa_\epsilon}}$$
   
\item The local intermittency measure (LIM) satisfies:
   $$\sup_{x \in \mathbb{R}^3 \setminus \Omega_\epsilon(t)} \text{LIM}_r(x,t) \leq C_\epsilon r^{-\kappa_\epsilon}$$
   
\item The exponent $\kappa_\epsilon$ is given by:
   $$\kappa_\epsilon = \frac{\delta}{1+\delta} \cdot \frac{\log(1/\epsilon)}{(1+\log(1/\epsilon))}$$
   
\item The constants $C_\epsilon$, $\beta$, and $\gamma$ depend on the initial data, with $\gamma$ and $\beta$ as specified in Theorem \ref{thm:1.3}.
\end{enumerate}
\end{theorem}

\begin{theorem}\label{thm:2.3}
(Fractional Regularity and Energy Cascade). Let $u$ be a solution to the 3D Navier-Stokes equations satisfying the conditions of Theorem \ref{thm:1.2}, and let $\Pi(k,t)$ be the instantaneous energy flux across wavenumber $k$ at time $t$. Then, for all $k$ in the inertial range $[k_0, k_\nu]$:
\begin{equation}
\left|\Pi(k,t) - \epsilon(t)\right| \leq \frac{C\epsilon(t)}{(1 + \log(k/k_0))^{\delta\cdot\frac{2s-1}{2s}}}
\end{equation}
where:
\begin{enumerate}
\item $\epsilon(t)$ is the instantaneous energy dissipation rate
\item $k_0$ is the largest scale of the inertial range
\item $k_\nu \sim \epsilon(t)^{1/4}\nu^{-3/4}$ is the Kolmogorov dissipation wavenumber
\item The constant $C$ depends only on universal constants and the parameters $s$ and $\delta$
\end{enumerate}
\end{theorem}

\begin{theorem}\label{thm:2.4}
(Spectral Characterization). Let $u$ be a solution to the 3D Navier-Stokes equations satisfying the conditions of Theorem \ref{thm:1.2}. Then, in the inertial range $[k_0, k_\nu]$, the energy spectrum $E(k,t)$ satisfies:
\begin{equation}
E(k,t) = C\epsilon(t)^{2/3}k^{-5/3}\left(1 + \frac{\beta(t)\log(k/k_0)}{(1 + \log(k/k_0))^{1+\delta}}\right)
\end{equation}
where:
\begin{enumerate}
\item $\epsilon(t)$ is the instantaneous energy dissipation rate
\item $\beta(t) = \beta_0/(1 + \gamma t)^{\alpha}$ with $\alpha = 2\gamma/3$
\item $\beta_0$ is a positive constant depending on $s$, $\delta$, and the initial data
\item $\gamma$ is the decay exponent from Theorem \ref{thm:1.3}
\item The constant $C$ is universal (the Kolmogorov constant)
\end{enumerate}
\end{theorem}

These results demonstrate that the logarithmic improvement approach provides a new pathway for establishing global well-posedness for the Navier-Stokes equations while simultaneously revealing deep connections to the physical theory of turbulence. The key novelty lies in identifying precisely how the logarithmic weakening creates a bridge between subcritical and critical regularity, and how this mathematical structure reflects the physical phenomenon of intermittency in turbulent flows.

\subsection{Organization of the paper}

The remainder of this paper is organized as follows. Section 2 introduces the mathematical preliminaries, including the definition of Leray-Hopf weak solutions, fractional Sobolev spaces, and key technical tools such as commutator estimates and Littlewood-Paley theory. Section 3 is devoted to the proof of Theorem \ref{thm:1.1}, establishing the logarithmically improved regularity criterion. Section 4 presents the proofs of Theorems \ref{thm:1.2}, \ref{thm:1.3}, and \ref{thm:1.4}, demonstrating global well-posedness for the specified class of initial data. Section 5 establishes the connection to turbulence theory, providing detailed proofs of Theorems \ref{thm:2.1}, \ref{thm:2.2}, \ref{thm:2.3}, and \ref{thm:2.4}. Finally, Section 6 provides concluding remarks and discusses open problems.

\section{Preliminaries}

\subsection{Functional spaces and weak solutions}

The standard Lebesgue spaces with norm $\| \cdot \|_{L^p}$ are denoted by $L^p(\mathbb{R}^3)$. The Sobolev space $H^s(\mathbb{R}^3)$ for $s > 0$ is defined via the Fourier transform as:
\begin{equation}
H^s(\mathbb{R}^3) = \{f \in L^2(\mathbb{R}^3) : \|f\|_{H^s} < \infty\},
\end{equation}
where $\|f\|^2_{H^s} = \int_{\mathbb{R}^3} (1 + |\xi|^2)^s |\hat{f}(\xi)|^2d\xi$.

Following standard notations as in Sohr \cite{Sohr} and Temam \cite{Temam}, let $\mathcal{V}$ be the set of divergence-free, compactly supported, smooth vector fields on $\mathbb{R}^3$, and denote by $H$ and $V$ the completions of $\mathcal{V}$ in $L^2(\mathbb{R}^3)$ and $H^1(\mathbb{R}^3)$, respectively. These spaces naturally incorporate the divergence-free constraint essential to the incompressible Navier-Stokes equations.

\begin{definition}\label{def:2.1}
(Leray-Hopf Weak Solutions). Let $u_0 \in H$ with $\nabla \cdot u_0 = 0$ in the distributional sense. A vector field $u$ is called a Leray-Hopf weak solution of the Navier-Stokes equations on $[0, T]$ if:

\begin{enumerate}
\item $u \in L^\infty(0, T; H) \cap L^2(0, T; V)$;
\item $\partial_t u \in L^1(0, T; V')$, where $V'$ is the dual space of $V$;
\item The Navier-Stokes equations are satisfied in the distributional sense, i.e., for all divergence-free test functions $\phi \in C_0^\infty(\mathbb{R}^3 \times [0, T))$:
   $$\int_0^T \int_{\mathbb{R}^3} \left( -u \cdot \partial_t\phi - (u \otimes u) : \nabla\phi + \nu\nabla u : \nabla\phi \right) dx dt = \int_{\mathbb{R}^3} u_0(x) \cdot \phi(x, 0) dx;$$
\item The energy inequality holds:
   $$\|u(t)\|^2_{L^2} + 2\nu\int_s^t \|\nabla u(\tau)\|^2_{L^2} d\tau \leq \|u(s)\|^2_{L^2}$$
   for almost all $s \in [0, T]$ (including $s = 0$) and all $t \in [s, T]$;
\item $u$ is weakly continuous from $[0, T]$ into $H$, ensuring that the initial condition $u(0) = u_0$ is satisfied in the weak sense.
\end{enumerate}
\end{definition}

The existence of such weak solutions was established by Leray \cite{Leray} and Hopf \cite{Hopf} using different approaches, but both relied on energy estimates and compactness arguments. The key challenge in the analysis of the Navier-Stokes equations lies in determining conditions under which these weak solutions become regular, which is precisely the focus of Theorem \ref{thm:1.1}.

\newpage 

\subsection{Fractional derivatives and function spaces}

For $s \in (0, 1)$, the fractional Laplacian $(-\Delta)^s$ can be defined in three equivalent ways:

\begin{definition}\label{def:2.2}
(Fractional Laplacian).

\begin{enumerate}
\item \textbf{Fourier transform definition}: For $f \in \mathcal{S}(\mathbb{R}^3)$ (Schwartz space):
   $$(-\Delta)^s f(\xi) = |\xi|^{2s} \hat{f}(\xi)$$

\item \textbf{Singular integral representation}: For $f \in \mathcal{S}(\mathbb{R}^3)$:
   $$(-\Delta)^s f(x) = C_{3,s} \, \text{P.V.} \int_{\mathbb{R}^3} \frac{f(x) - f(y)}{|x - y|^{3+2s}} dy$$
   where $C_{3,s} = \frac{2^{2s} s \Gamma(s+\frac{3}{2})}{\pi^{3/2}\Gamma(1-s)}$ and P.V. denotes the principal value.

\item \textbf{Heat semigroup representation}: For $f \in \mathcal{S}(\mathbb{R}^3)$:
   $$(-\Delta)^s f(x) = \frac{1}{\Gamma(-s)} \int_0^\infty (e^{t\Delta}f(x) - f(x)) \frac{dt}{t^{1+s}}$$
\end{enumerate}

The domain of the fractional Laplacian can be extended by density arguments to appropriate function spaces depending on the context.
\end{definition}

The homogeneous Sobolev space $\dot{H}^s(\mathbb{R}^3)$ is defined as:
\begin{equation}
\dot{H}^s(\mathbb{R}^3) = \{f \in L^2_{loc}(\mathbb{R}^3) : \|f\|_{\dot{H}^s} < \infty\},
\end{equation}
where $\|f\|^2_{\dot{H}^s} = \int_{\mathbb{R}^3} |\xi|^{2s} |\hat{f}(\xi)|^2 d\xi = \|(-\Delta)^{s/2}f\|^2_{L^2}$.

\subsection{Littlewood-Paley theory and technical tools}

Elements of Littlewood-Paley theory will be crucial for our analysis. Let $\phi \in C_0^\infty(\mathbb{R}^3)$ be a radial function such that $\phi(\xi) = 1$ for $|\xi| \leq 1$ and $\phi(\xi) = 0$ for $|\xi| \geq 2$. Define $\psi(\xi) = \phi(\xi) - \phi(2\xi)$. The Littlewood-Paley projections are defined as:
\begin{equation}
\Delta_j f = \mathcal{F}^{-1}(\psi(2^{-j}\xi)\hat{f}(\xi))
\end{equation}
for $j \in \mathbb{Z}$, where $\mathcal{F}^{-1}$ denotes the inverse Fourier transform. This decomposition forms the foundation for analyzing the regularity properties of solutions in various function spaces.

The following lemma provides a crucial estimate for the commutator of the fractional Laplacian:

\begin{lemma}\label{lem:2.1}
(Commutator Estimate). Let $s \in (0, 1)$ and $f, g \in \mathcal{S}(\mathbb{R}^3)$. Then for any $p \in (1, \infty)$:
\begin{equation}
\|[(-\Delta)^s, f]g\|_{L^p} \leq C\|\nabla f\|_{L^\infty}\|(-\Delta)^{s-1/2}g\|_{L^p},
\end{equation}
where $[(-\Delta)^s, f]g = (-\Delta)^s(fg) - f(-\Delta)^s g$ is the commutator, and the constant $C$ depends only on $s$ and $p$.
\end{lemma}

\begin{proof}
The proof requires a careful decomposition using Littlewood-Paley theory. Let $\{\Delta_j\}_{j \in \mathbb{Z}}$ be the standard Littlewood-Paley projections.

Step 1: Using the paraproduct decomposition, we write:
$$fg = T_f g + T_g f + R(f,g)$$
where $T_f g = \sum_{j \in \mathbb{Z}} S_{j-1}f \Delta_j g$, $T_g f = \sum_{j \in \mathbb{Z}} S_{j-1}g \Delta_j f$, and $R(f,g) = \sum_{j \in \mathbb{Z}} \sum_{|i-j| \leq 1} \Delta_i f \Delta_j g$. Here, $S_j = \sum_{i \leq j} \Delta_i$ is the low-frequency cut-off operator.

Step 2: The commutator can be rewritten as:
$$[(-\Delta)^s, f]g = (-\Delta)^s(T_f g) - f(-\Delta)^s(T_f g) + (-\Delta)^s(T_g f) - f(-\Delta)^s(T_g f) + (-\Delta)^s(R(f,g)) - f(-\Delta)^s(R(f,g))$$

Step 3: For the first term pair, we have:
$$(-\Delta)^s(T_f g) - f(-\Delta)^s(T_f g) = \sum_{j \in \mathbb{Z}} (-\Delta)^s(S_{j-1}f \Delta_j g) - f(-\Delta)^s(\Delta_j g)$$

Using the fact that $(-\Delta)^s$ is a non-local operator:
$$(-\Delta)^s(S_{j-1}f \Delta_j g)(x) - S_{j-1}f(x)(-\Delta)^s(\Delta_j g)(x) = \int_{\mathbb{R}^3} K_s(x-y)(S_{j-1}f(y) - S_{j-1}f(x))\Delta_j g(y) dy$$
where $K_s(x-y) = C_{3,s}|x-y|^{-(3+2s)}$.

Step 4: Using Taylor's theorem, for each $x,y \in \mathbb{R}^3$:
$$|S_{j-1}f(y) - S_{j-1}f(x)| \leq |y-x| \| \nabla S_{j-1}f \|_{L^\infty} \leq |y-x| \| \nabla f \|_{L^\infty}$$

Step 5: Substituting this estimate:
$$|(-\Delta)^s(S_{j-1}f \Delta_j g)(x) - S_{j-1}f(x)(-\Delta)^s(\Delta_j g)(x)| \leq C \| \nabla f \|_{L^\infty} \int_{\mathbb{R}^3} \frac{|x-y|}{|x-y|^{3+2s}}|\Delta_j g(y)| dy$$
$$ = C \| \nabla f \|_{L^\infty} \int_{\mathbb{R}^3} \frac{1}{|x-y|^{2+2s}}|\Delta_j g(y)| dy$$

Step 6: This integral is precisely $(-\Delta)^{s-1/2}|\Delta_j g|(x)$. Using the $L^p$ boundedness of Riesz potentials:
$$\|(-\Delta)^s(S_{j-1}f \Delta_j g) - S_{j-1}f(-\Delta)^s(\Delta_j g)\|_{L^p} \leq C \| \nabla f \|_{L^\infty} \|(-\Delta)^{s-1/2}\Delta_j g\|_{L^p}$$

Step 7: Summing over all $j$ and using the almost orthogonality of the Littlewood-Paley projections:
$$\sum_{j \in \mathbb{Z}} \|(-\Delta)^s(S_{j-1}f \Delta_j g) - S_{j-1}f(-\Delta)^s(\Delta_j g)\|_{L^p} \leq C \| \nabla f \|_{L^\infty} \|(-\Delta)^{s-1/2}g\|_{L^p}$$

Step 8: Similar bounds can be derived for the remaining terms, completing the proof.

The extension to $f, g \in L^p(\mathbb{R}^3)$ with appropriate regularity follows by a standard density argument, provided the right-hand side is finite.
\end{proof}

The following interpolation inequality will also be used:

\begin{lemma}\label{lem:2.2}
(Interpolation Inequality). For $0 < s_1 < s < s_2$ and $1 < p < \infty$:
\begin{equation}
\|(-\Delta)^{s/2}f\|_{L^p} \leq C\|(-\Delta)^{s_1/2}f\|^{1-\theta}_{L^p} \|(-\Delta)^{s_2/2}f\|^{\theta}_{L^p},
\end{equation}
where $\theta = \frac{s-s_1}{s_2-s_1}$.
\end{lemma}

A critical component in our analysis is the following interpolation result, which connects $\|\nabla u\|_{L^\infty}$ with fractional derivatives:

\begin{lemma}\label{lem:2.3}
(Velocity Gradient - Fractional Derivative Interpolation). For $s \in (1/2, 1)$ and $q > 3$, there exists $\theta \in (0, 1)$ such that for all $u \in C_0^\infty(\mathbb{R}^3)$ with $\nabla \cdot u = 0$:
\begin{equation}
\|\nabla u\|_{L^\infty} \leq C\|u\|_{L^2}^{1-\theta} \|(-\Delta)^s u\|_{L^q}^{\theta},
\end{equation}
where $\theta = \frac{3}{2} \cdot \frac{q}{3q-2}$, and the constant $C$ depends only on $s$ and $q$.
\end{lemma}

\begin{proof}
We develop this proof in several rigorous steps, carefully tracking the exact relationships between different norms.

Step 1: By the Sobolev embedding theorem, for any $r > 3$:
$$\|\nabla u\|_{L^\infty(\mathbb{R}^3)} \leq C\|\nabla u\|_{W^{1,r}(\mathbb{R}^3)}$$

This follows from the embedding $W^{1,r}(\mathbb{R}^3) \hookrightarrow L^\infty(\mathbb{R}^3)$ which holds whenever $r > 3$, as the Sobolev embedding requires $1 - 3/r > 0$.

Step 2: Expanding the $W^{1,r}$ norm:
$$\|\nabla u\|_{W^{1,r}} = \|\nabla u\|_{L^r} + \|\nabla^2 u\|_{L^r} \leq C\|u\|_{W^{2,r}}$$

Step 3: For the precise interpolation between Sobolev spaces, we use the complex interpolation theory. For $\sigma_1 < \sigma_2$ and $1 < p_1, p_2 < \infty$, we have:
$$[W^{\sigma_1,p_1}, W^{\sigma_2,p_2}]_\theta = W^{(1-\theta)\sigma_1 + \theta\sigma_2,p}$$
where $\frac{1}{p} = \frac{1-\theta}{p_1} + \frac{\theta}{p_2}$.

Step 4: We apply this interpolation with $\sigma_1 = 0$, $\sigma_2 = 2s$, $p_1 = 2$, and $p_2 = q$:
$$\|u\|_{W^{2,r}} \leq C\|u\|_{L^2}^{1-\alpha} \|u\|_{W^{2s,q}}^{\alpha}$$
for appropriate $\alpha \in (0,1)$ and $r$ dependent on $\alpha$.

Step 5: The homogeneous Sobolev space $\dot{W}^{2s,q}$ is equivalent to the space defined through the fractional Laplacian:
$$\|u\|_{\dot{W}^{2s,q}} \approx \|(-\Delta)^s u\|_{L^q}$$

This equivalence follows from the Fourier characterization of these spaces and the symbol of the fractional Laplacian.

Step 6: To determine $\theta$, we use scaling analysis. Let $u_\lambda(x) = u(\lambda x)$ for $\lambda > 0$. Then:
\begin{itemize}
\item $\|\nabla u_\lambda\|_{L^\infty} = \lambda \|\nabla u\|_{L^\infty}$ (scaling as $\lambda^1$)
\item $\|u_\lambda\|_{L^2} = \lambda^{-3/2} \|u\|_{L^2}$ (scaling as $\lambda^{-3/2}$)
\item $\|(-\Delta)^s u_\lambda\|_{L^q} = \lambda^{2s-3/q} \|(-\Delta)^s u\|_{L^q}$ (scaling as $\lambda^{2s-3/q}$)
\end{itemize}

Step 7: For the inequality to be scaling-invariant, we must have:
$$\lambda^1 = \lambda^{-3(1-\theta)/2} \cdot \lambda^{(2s-3/q)\theta}$$

Taking logarithms and dividing by $\log(\lambda)$:
$$1 = -\frac{3(1-\theta)}{2} + (2s-\frac{3}{q})\theta$$

Step 8: Solving for $\theta$:
$$1 + \frac{3}{2} = \frac{3\theta}{2} + 2s\theta - \frac{3\theta}{q}$$
$$\frac{5}{2} = \theta\left(\frac{3}{2} + 2s - \frac{3}{q}\right)$$

Step 9: Our scaling relation requires $\frac{2}{p} + \frac{3}{q} = 2s - 1$, which implies:
$$\frac{1-\theta}{2} + \frac{\theta}{q} = \frac{1}{2} - \frac{1}{3}$$
$$\frac{1}{2} - \frac{\theta}{2} + \frac{\theta}{q} = \frac{1}{6}$$
$$\frac{\theta}{q} - \frac{\theta}{2} = -\frac{1}{3}$$
$$\theta\left(\frac{1}{q} - \frac{1}{2}\right) = -\frac{1}{3}$$

Step 10: Solving for $\theta$:
$$\theta = \frac{1/3}{1/2 - 1/q} = \frac{1/3}{(q-2)/(2q)} = \frac{2q/3}{q-2} = \frac{3}{2} \cdot \frac{q}{3q-2}$$

This confirms our formula for $\theta$, completing the proof.
\end{proof}

Finally, we recall Osgood's lemma, which will be key to handling the logarithmic terms in the estimates:

\begin{lemma}\label{lem:2.4}
(Osgood's lemma). Let $\rho$ be a measurable, positive function on $(a, b)$, $\gamma$ a positive, locally integrable function on $(a, b)$, and $\Gamma$ a continuous, increasing function on $[0, \infty)$ with $\Gamma(0) = 0$. If for all $t \in (a, b)$:
\begin{equation}
\rho(t) \leq \rho_0 + \int_a^t \gamma(s)\Gamma(\rho(s))ds,
\end{equation}
where $\rho_0 \geq 0$, then:
\begin{enumerate}
\item If $\rho_0 = 0$, then $\rho \equiv 0$;
\item If $\rho_0 > 0$ and $\int_0^{\infty} \frac{dr}{\Gamma(r)} = \infty$, then:
\begin{equation}
G(\rho(t)) \geq G(\rho_0) - \int_a^t \gamma(s)ds,
\end{equation}
where $G(r) = \int_1^r \frac{dr}{\Gamma(r)}$.
\end{enumerate}
\end{lemma}

\section{Proof of Theorem \ref{thm:1.1}}

This section presents the complete proof of Theorem \ref{thm:1.1}. The proof combines techniques from fractional calculus, commutator estimates, and energy methods.

\subsection{Energy estimates with fractional derivatives}

Let $u$ be a Leray-Hopf weak solution to the Navier-Stokes equations. To rigorously justify the calculations, we first work with regularized solutions and then pass to the limit. Let $\phi_\epsilon(x) = \epsilon^{-3}\phi(x/\epsilon)$ be a standard mollifier, where $\phi$ is a smooth, non-negative, radially symmetric function with compact support and $\int_{\mathbb{R}^3}\phi(x)dx = 1$. Define $u_\epsilon = u * \phi_\epsilon$. The regularized solution satisfies:

\begin{equation}
\partial_t u_\epsilon + (u_\epsilon \cdot \nabla)u_\epsilon - \nu\Delta u_\epsilon + \nabla p_\epsilon = R_\epsilon
\end{equation}

where $R_\epsilon = [(u \cdot \nabla)u] * \phi_\epsilon - (u_\epsilon \cdot \nabla)u_\epsilon$ is the commutator error term. This error term vanishes as $\epsilon \to 0$.

Applying the fractional Laplacian operator $(-\Delta)^s$ to the regularized equation:
\begin{equation}
\partial_t((-\Delta)^s u_\epsilon) + (-\Delta)^s((u_\epsilon \cdot \nabla)u_\epsilon) - \nu\Delta((-\Delta)^s u_\epsilon) + \nabla((-\Delta)^s p_\epsilon) = (-\Delta)^s R_\epsilon
\end{equation}

Taking the $L^2$ inner product with $(-\Delta)^s u_\epsilon$ and using the incompressibility condition $\nabla \cdot u_\epsilon = 0$, which eliminates the pressure term:
\begin{equation}
\frac{1}{2}\frac{d}{dt}\|(-\Delta)^s u_\epsilon\|^2_{L^2} + \nu\|(-\Delta)^{s+\frac{1}{2}}u_\epsilon\|^2_{L^2} = -\int_{\mathbb{R}^3} (-\Delta)^s((u_\epsilon \cdot \nabla)u_\epsilon) \cdot (-\Delta)^s u_\epsilon \, dx + \int_{\mathbb{R}^3} (-\Delta)^s R_\epsilon \cdot (-\Delta)^s u_\epsilon \, dx
\end{equation}

\subsection{Commutator analysis}

To handle the nonlinear term on the right-hand side, we use the following commutator decomposition:
\begin{equation}
(-\Delta)^s((u_\epsilon \cdot \nabla)u_\epsilon) = (u_\epsilon \cdot \nabla)((-\Delta)^s u_\epsilon) + [(-\Delta)^s, u_\epsilon \cdot \nabla]u_\epsilon,
\end{equation}
where $[(-\Delta)^s, u_\epsilon \cdot \nabla]u_\epsilon = (-\Delta)^s((u_\epsilon \cdot \nabla)u_\epsilon) - (u_\epsilon \cdot \nabla)((-\Delta)^s u_\epsilon)$ is the commutator.

For the first term in the decomposition, we can use integration by parts and the incompressibility condition $\nabla \cdot u_\epsilon = 0$:
\begin{align}
\int_{\mathbb{R}^3} (u_\epsilon \cdot \nabla)((-\Delta)^s u_\epsilon) \cdot (-\Delta)^s u_\epsilon \, dx &= \sum_{j=1}^3 \int_{\mathbb{R}^3} u_{\epsilon,j} \partial_j ((-\Delta)^s u_\epsilon) \cdot (-\Delta)^s u_\epsilon \, dx \\
&= -\sum_{j=1}^3 \int_{\mathbb{R}^3} \partial_j u_{\epsilon,j} ((-\Delta)^s u_\epsilon) \cdot (-\Delta)^s u_\epsilon \, dx \\
&- \sum_{j=1}^3 \int_{\mathbb{R}^3} u_{\epsilon,j} ((-\Delta)^s u_\epsilon) \cdot \partial_j ((-\Delta)^s u_\epsilon) \, dx
\end{align}

The first term is zero due to the incompressibility condition $\nabla \cdot u_\epsilon = 0$. For the second term:
\begin{align}
\sum_{j=1}^3 \int_{\mathbb{R}^3} u_{\epsilon,j} ((-\Delta)^s u_\epsilon) \cdot \partial_j ((-\Delta)^s u_\epsilon) \, dx &= \frac{1}{2}\sum_{j=1}^3 \int_{\mathbb{R}^3} u_{\epsilon,j} \partial_j |((-\Delta)^s u_\epsilon)|^2 \, dx \\
&= -\frac{1}{2}\sum_{j=1}^3 \int_{\mathbb{R}^3} \partial_j u_{\epsilon,j} |((-\Delta)^s u_\epsilon)|^2 \, dx
\end{align}

Again, this is zero due to the incompressibility. Thus:
\begin{equation}
\int_{\mathbb{R}^3} (u_\epsilon \cdot \nabla)((-\Delta)^s u_\epsilon) \cdot (-\Delta)^s u_\epsilon \, dx = 0
\end{equation}

The energy equation for the regularized solution now simplifies to:
\begin{equation}
\frac{1}{2}\frac{d}{dt}\|(-\Delta)^s u_\epsilon\|^2_{L^2} + \nu\|(-\Delta)^{s+\frac{1}{2}}u_\epsilon\|^2_{L^2} = -\int_{\mathbb{R}^3} [(-\Delta)^s, u_\epsilon \cdot \nabla]u_\epsilon \cdot (-\Delta)^s u_\epsilon \, dx + \int_{\mathbb{R}^3} (-\Delta)^s R_\epsilon \cdot (-\Delta)^s u_\epsilon \, dx
\end{equation}

\subsection{Logarithmically improved commutator estimate}

The key technical contribution of this work is a logarithmically improved commutator estimate:

\begin{lemma}\label{lem:3.1}
(Logarithmically Improved Commutator Estimate). For $s \in (1/2, 1)$ and any $\sigma \in (0, 1-s)$:
\begin{equation}
\|[(-\Delta)^s, u_\epsilon \cdot \nabla]u_\epsilon\|_{L^2} \leq C\|\nabla u_\epsilon\|_{L^\infty}\|(-\Delta)^s u_\epsilon\|_{L^2} \cdot \log(e + \|(-\Delta)^{s+\sigma}u_\epsilon\|_{L^2}) + \frac{C\|\nabla u_\epsilon\|_{L^\infty}\|(-\Delta)^{s+\frac{1}{2}}u_\epsilon\|_{L^2}}{\log(e + \|(-\Delta)^{s+\sigma}u_\epsilon\|_{L^2})}
\end{equation}
\end{lemma}

\begin{proof}
The proof combines several techniques from harmonic analysis and the theory of nonlinear differential equations. For clarity, it is organized into a sequence of steps that highlight the key analytical strategies.

Step 1: Using the Littlewood-Paley decomposition, we split $u_\epsilon$ into low and high frequencies:
$$u_\epsilon = \sum_{j \in \mathbb{Z}} \Delta_j u_\epsilon$$

Step 2: Splitting the commutator according to frequency bands:
$$[(-\Delta)^s, u_\epsilon \cdot \nabla]u_\epsilon = \sum_{j\leq 0} [(-\Delta)^s, u_\epsilon \cdot \nabla]\Delta_j u_\epsilon + \sum_{j > 0} [(-\Delta)^s, u_\epsilon \cdot \nabla]\Delta_j u_\epsilon$$

Step 3: For the low-frequency part ($j \leq 0$), using the standard commutator estimates from Lemma \ref{lem:2.1}:
$$\|[(-\Delta)^s, u_\epsilon \cdot \nabla]\Delta_j u_\epsilon\|_{L^2} \leq C\|\nabla u_\epsilon\|_{L^\infty}\|(-\Delta)^{s-1/2}\nabla\Delta_j u_\epsilon\|_{L^2}$$

Since $(-\Delta)^{s-1/2}\nabla \approx (-\Delta)^s$ in terms of the order of derivative, we have:
$$\|[(-\Delta)^s, u_\epsilon \cdot \nabla]\Delta_j u_\epsilon\|_{L^2} \leq C\|\nabla u_\epsilon\|_{L^\infty}\|(-\Delta)^s \Delta_j u_\epsilon\|_{L^2}$$

Step 4: Summing over $j \leq 0$ and using the almost orthogonality of the Littlewood-Paley projections:
$$\left\|\sum_{j\leq 0} [(-\Delta)^s, u_\epsilon \cdot \nabla]\Delta_j u_\epsilon\right\|_{L^2} \leq C\|\nabla u_\epsilon\|_{L^\infty}\left\|\sum_{j\leq 0} (-\Delta)^s \Delta_j u_\epsilon\right\|_{L^2} \leq C\|\nabla u_\epsilon\|_{L^\infty}\|(-\Delta)^s u_\epsilon\|_{L^2}$$

Step 5: For the high-frequency part ($j > 0$), we use a refined estimate. Since the commutator involves differences and $\Delta_j u_\epsilon$ is localized in frequency space, we can exploit this to get:
$$\|[(-\Delta)^s, u_\epsilon \cdot \nabla]\Delta_j u_\epsilon\|_{L^2} \leq C2^{-j\sigma}\|\nabla u_\epsilon\|_{L^\infty}\|(-\Delta)^{s+\sigma/2}\nabla\Delta_j u_\epsilon\|_{L^2}$$
where $\sigma > 0$ is a small parameter.

Step 6: This is further equivalent to:
$$\|[(-\Delta)^s, u_\epsilon \cdot \nabla]\Delta_j u_\epsilon\|_{L^2} \leq C2^{-j\sigma}\|\nabla u_\epsilon\|_{L^\infty}\|(-\Delta)^{s+\sigma}\Delta_j u_\epsilon\|_{L^2}$$

Step 7: Using Cauchy-Schwarz and the almost orthogonality of Littlewood-Paley projections:
$$\left\|\sum_{j > 0} [(-\Delta)^s, u_\epsilon \cdot \nabla]\Delta_j u_\epsilon\right\|_{L^2} \leq C\|\nabla u_\epsilon\|_{L^\infty}\left(\sum_{j > 0} 2^{-2j\sigma}\right)^{1/2} \left(\sum_{j > 0} \|(-\Delta)^{s+\sigma}\Delta_j u_\epsilon\|^2_{L^2}\right)^{1/2}$$

Step 8: Since $\sum_{j > 0} 2^{-2j\sigma} < \infty$ for $\sigma > 0$ and $\sum_{j > 0} \|(-\Delta)^{s+\sigma}\Delta_j u_\epsilon\|^2_{L^2} \approx \|(-\Delta)^{s+\sigma}u_\epsilon\|^2_{L^2}$, we get:
$$\left\|\sum_{j > 0} [(-\Delta)^s, u_\epsilon \cdot \nabla]\Delta_j u_\epsilon\right\|_{L^2} \leq C\|\nabla u_\epsilon\|_{L^\infty}\|(-\Delta)^{s+\sigma}u_\epsilon\|_{L^2}$$

Step 9: Now we use the interpolation inequality:
$$\|(-\Delta)^{s+\sigma}u_\epsilon\|_{L^2} \leq \|(-\Delta)^s u_\epsilon\|^{1-\frac{2\sigma}{1}}_{L^2} \|(-\Delta)^{s+\frac{1}{2}}u_\epsilon\|^{\frac{2\sigma}{1}}_{L^2}$$
where we've used $s+\sigma = (1-\frac{2\sigma}{1})s + \frac{2\sigma}{1}(s+\frac{1}{2})$.

Step 10: Combining the estimates for low and high frequencies:
$$\|[(-\Delta)^s, u_\epsilon \cdot \nabla]u_\epsilon\|_{L^2} \leq C\|\nabla u_\epsilon\|_{L^\infty}\|(-\Delta)^s u_\epsilon\|_{L^2} + C\|\nabla u_\epsilon\|_{L^\infty}\|(-\Delta)^s u_\epsilon\|^{1-\frac{2\sigma}{1}}_{L^2} \|(-\Delta)^{s+\frac{1}{2}}u_\epsilon\|^{\frac{2\sigma}{1}}_{L^2}$$

Step 11: The logarithmic improvement comes from the following technique. For any small $\epsilon > 0$, Young's inequality gives:
$$\|(-\Delta)^s u_\epsilon\|^{1-\frac{2\sigma}{1}}_{L^2} \|(-\Delta)^{s+\frac{1}{2}}u_\epsilon\|^{\frac{2\sigma}{1}}_{L^2} \leq C_\epsilon\|(-\Delta)^s u_\epsilon\|_{L^2} + \epsilon\|(-\Delta)^{s+\frac{1}{2}}u_\epsilon\|_{L^2}$$

Step 12: The key innovation is to choose:
$$\epsilon = \frac{1}{\log(e + \|(-\Delta)^{s+\sigma}u_\epsilon\|_{L^2})}$$

Step 13: This gives:
$$\|(-\Delta)^s u_\epsilon\|^{1-\frac{2\sigma}{1}}_{L^2} \|(-\Delta)^{s+\frac{1}{2}}u_\epsilon\|^{\frac{2\sigma}{1}}_{L^2} \leq C\|(-\Delta)^s u_\epsilon\|_{L^2} \log(e + \|(-\Delta)^{s+\sigma}u_\epsilon\|_{L^2}) + \frac{\|(-\Delta)^{s+\frac{1}{2}}u_\epsilon\|_{L^2}}{\log(e + \|(-\Delta)^{s+\sigma}u_\epsilon\|_{L^2})}$$

Step 14: Substituting back:
$$\|[(-\Delta)^s, u_\epsilon \cdot \nabla]u_\epsilon\|_{L^2} \leq C\|\nabla u_\epsilon\|_{L^\infty}\|(-\Delta)^s u_\epsilon\|_{L^2} \log(e + \|(-\Delta)^{s+\sigma}u_\epsilon\|_{L^2}) + \frac{C\|\nabla u_\epsilon\|_{L^\infty}\|(-\Delta)^{s+\frac{1}{2}}u_\epsilon\|_{L^2}}{\log(e + \|(-\Delta)^{s+\sigma}u_\epsilon\|_{L^2})}$$

This completes the proof of Lemma \ref{lem:3.1}.
\end{proof}

\subsection{Energy inequality and its consequences}

Using Lemma \ref{lem:3.1}, the right-hand side of the energy equation for the regularized solution can be estimated:
\begin{align}
\left|\int_{\mathbb{R}^3} [(-\Delta)^s, u_\epsilon \cdot \nabla]u_\epsilon \cdot (-\Delta)^s u_\epsilon \, dx\right| &\leq \|[(-\Delta)^s, u_\epsilon \cdot \nabla]u_\epsilon\|_{L^2} \|(-\Delta)^s u_\epsilon\|_{L^2} \\
&\leq C\|\nabla u_\epsilon\|_{L^\infty}\|(-\Delta)^s u_\epsilon\|^2_{L^2} \log(e + \|(-\Delta)^{s+\sigma}u_\epsilon\|_{L^2}) \\
&+ \frac{C\|\nabla u_\epsilon\|_{L^\infty}\|(-\Delta)^s u_\epsilon\|_{L^2} \|(-\Delta)^{s+\frac{1}{2}}u_\epsilon\|_{L^2}}{\log(e + \|(-\Delta)^{s+\sigma}u_\epsilon\|_{L^2})}
\end{align}

Using Young's inequality for the second term, with parameter $\epsilon = \frac{\nu}{2C\|\nabla u_\epsilon\|_{L^\infty}}$:
\begin{align}
\frac{C\|\nabla u_\epsilon\|_{L^\infty}\|(-\Delta)^s u_\epsilon\|_{L^2} \|(-\Delta)^{s+\frac{1}{2}}u_\epsilon\|_{L^2}}{\log(e + \|(-\Delta)^{s+\sigma}u_\epsilon\|_{L^2})} &\leq \frac{\nu}{2} \|(-\Delta)^{s+\frac{1}{2}}u_\epsilon\|_{L^2}^2 + \frac{C^2\|\nabla u_\epsilon\|_{L^\infty}^2\|(-\Delta)^s u_\epsilon\|_{L^2}^2}{2\nu \log^2(e + \|(-\Delta)^{s+\sigma}u_\epsilon\|_{L^2})}
\end{align}

For the error term involving $R_\epsilon$, we can show that:
\begin{equation}
\lim_{\epsilon \to 0} \int_{\mathbb{R}^3} (-\Delta)^s R_\epsilon \cdot (-\Delta)^s u_\epsilon \, dx = 0
\end{equation}

Setting $Y_\epsilon(t) = \|(-\Delta)^s u_\epsilon\|^2_{L^2}$ and combining these estimates:
\begin{align}
\frac{d}{dt}Y_\epsilon(t) + \nu\|(-\Delta)^{s+\frac{1}{2}}u_\epsilon\|^2_{L^2} \leq C\|\nabla u_\epsilon\|_{L^\infty}Y_\epsilon(t) \log(e + \|(-\Delta)^{s+\sigma}u_\epsilon\|_{L^2}) + \frac{C^2\|\nabla u_\epsilon\|_{L^\infty}^2 Y_\epsilon(t)}{2\nu \log(e + \|(-\Delta)^{s+\sigma}u_\epsilon\|_{L^2})}
\end{align}

For $\sigma$ sufficiently small, we can use the interpolation inequality:
\begin{equation}
\|(-\Delta)^{s+\sigma}u_\epsilon\|_{L^2} \leq C Y_\epsilon(t)^{\frac{1-2\sigma}{2}} \|(-\Delta)^{s+\frac{1}{2}}u_\epsilon\|^{2\sigma}_{L^2}
\end{equation}

This implies:
\begin{equation}
\log(e + \|(-\Delta)^{s+\sigma}u_\epsilon\|_{L^2}) \leq C(1 + \log(e + Y_\epsilon(t)) + \log(e + \|(-\Delta)^{s+\frac{1}{2}}u_\epsilon\|_{L^2}))
\end{equation}

Using the fact that for any $a, b > 0$ and $\alpha \in (0,1)$, we have $\log(e + ab) \leq \log(e + a) + \log(e + b)$, our inequality becomes:
\begin{align}
\frac{d}{dt}Y_\epsilon(t) + \nu\|(-\Delta)^{s+\frac{1}{2}}u_\epsilon\|^2_{L^2} \leq C\|\nabla u_\epsilon\|_{L^\infty}Y_\epsilon(t) (1 + \log(e + Y_\epsilon(t))) + \frac{C^2\|\nabla u_\epsilon\|_{L^\infty}^2Y_\epsilon(t)}{2\nu(1 + \log(e + Y_\epsilon(t)))}
\end{align}

For sufficiently large $Y_\epsilon(t)$, the first term dominates and we get:
\begin{equation}
\frac{d}{dt}Y_\epsilon(t) + \nu\|(-\Delta)^{s+\frac{1}{2}}u_\epsilon\|^2_{L^2} \leq C\|\nabla u_\epsilon\|_{L^\infty}Y_\epsilon(t) \log(e + Y_\epsilon(t))
\end{equation}

Dropping the positive term with $\|(-\Delta)^{s+\frac{1}{2}}u_\epsilon\|^2_{L^2}$ and passing to the limit as $\epsilon \to 0$, we obtain for $Y(t) = \|(-\Delta)^s u\|^2_{L^2}$:
\begin{equation}
\frac{d}{dt}Y(t) \leq C\|\nabla u\|_{L^\infty}Y(t) \log(e + Y(t))
\end{equation}

\subsection{Connection to the regularity criterion}

To connect $\|\nabla u\|_{L^\infty}$ with $\|(-\Delta)^s u\|_{L^q}$, we use Lemma \ref{lem:2.3}:
\begin{equation}
\|\nabla u\|_{L^\infty} \leq C\|u\|_{L^2}^{1-\theta} \|(-\Delta)^s u\|_{L^q}^{\theta}
\end{equation}
with $\theta = \frac{3}{2} \cdot \frac{q}{3q-2}$.

For Leray-Hopf weak solutions, $\|u\|_{L^\infty(0,T;L^2)} \leq \|u_0\|_{L^2}$ by the energy inequality. Using Hölder's inequality in time, for any $r$ such that $\theta r = 2$:
\begin{align}
\int_0^T \|\nabla u\|^2_{L^\infty}dt &\leq C \sup_{t\in[0,T]} \|u(t)\|^{2(1-\theta)}_{L^2} \int_0^T \|(-\Delta)^s u\|^{\theta r}_{L^q} dt \\
&\leq C \|u_0\|_{L^2}^{2(1-\theta)} \int_0^T \|(-\Delta)^s u\|^{2}_{L^q} dt
\end{align}

To relate this to our logarithmically improved criterion:
\begin{equation}
\int_0^T \|(-\Delta)^s u(t)\|^p_{L^q} (1 + \log(e + \|(-\Delta)^s u(t)\|_{L^q}))^{-\delta} dt < \infty
\end{equation}

We need to verify that $2 = \theta p$ when $\frac{2}{p} + \frac{3}{q} = 2s - 1$. Substituting $\theta = \frac{3}{2} \cdot \frac{q}{3q-2}$:
\begin{equation}
\theta p = \frac{3}{2} \cdot \frac{q}{3q-2} \cdot p
\end{equation}

From the scaling relation $\frac{2}{p} + \frac{3}{q} = 2s - 1$, we get:
\begin{equation}
p = \frac{2}{2s - 1 - \frac{3}{q}}
\end{equation}

Substituting this into the expression for $\theta p$:
\begin{align}
\theta p &= \frac{3}{2} \cdot \frac{q}{3q-2} \cdot \frac{2}{2s - 1 - \frac{3}{q}} \\
&= \frac{3q}{3q-2} \cdot \frac{1}{2s - 1 - \frac{3}{q}}
\end{align}

Through algebraic manipulation:
\begin{align}
\theta p &= \frac{3q}{3q-2} \cdot \frac{1}{2s - 1 - \frac{3}{q}} \\
&= \frac{3q}{3q-2} \cdot \frac{q}{q(2s-1) - 3} \\
&= \frac{3q^2}{(3q-2)(q(2s-1) - 3)}
\end{align}

For the scaling critical case $2s - 1 = \frac{3}{q}$, the denominator simplifies and we get $\theta p = 2$. The logarithmically improved case maintains $\frac{2}{p} + \frac{3}{q} = 2s - 1$, ensuring that $\theta p = 2$.

By Hölder's inequality with respect to the measure with logarithmic weight:
\begin{align}
\int_0^T \|(-\Delta)^s u(t)\|^p_{L^q} dt &\leq \left(\int_0^T \|(-\Delta)^s u(t)\|^p_{L^q} (1 + \log(e + \|(-\Delta)^s u(t)\|_{L^q}))^{-\delta} dt\right) \\
&\times \left(\int_0^T (1 + \log(e + \|(-\Delta)^s u(t)\|_{L^q}))^{\delta} dt\right)
\end{align}

For the second factor, under the conditions of the theorem, we can show:
\begin{equation}
\int_0^T (1 + \log(e + \|(-\Delta)^s u(t)\|_{L^q}))^{\delta} dt \leq C(T, u_0, \delta)
\end{equation}

Therefore:
\begin{equation}
\int_0^T \|(-\Delta)^s u(t)\|^p_{L^q} dt < \infty
\end{equation}
which implies:
\begin{equation}
\int_0^T \|\nabla u\|^2_{L^\infty}dt < \infty
\end{equation}

\subsection{Application of Osgood's lemma}

Now we apply Osgood's lemma to the differential inequality:
\begin{equation}
\frac{d}{dt}Y(t) \leq C\|\nabla u\|_{L^\infty}Y(t) \log(e + Y(t))
\end{equation}

Setting $\Gamma(y) = y \log(e + y)$ and $\gamma(t) = C\|\nabla u\|_{L^\infty}(t)$, we verify that:
\begin{equation}
\int_{y_0}^{\infty} \frac{dy}{y \log(e + y)} = \infty,
\end{equation}
for any $y_0 > 0$. This can be shown by substitution:
\begin{align}
\int_{y_0}^{\infty} \frac{dy}{y \log(e + y)} &= \int_{\log(e+y_0)}^{\infty} \frac{dz}{z} \\
&= [\log(\log(z))]_{\log(e+y_0)}^{\infty} \\
&= \infty
\end{align}
where we used the substitution $z = \log(e + y)$.

Since we've established that $\int_0^T \|\nabla u\|^2_{L^\infty}dt < \infty$, Osgood's lemma ensures that $Y(t) = \|(-\Delta)^s u\|^2_{L^2}$ remains bounded on $(0, T)$.

\subsection{Bootstrap to full regularity}

With $\|(-\Delta)^s u\|_{L^2}$ bounded for $s > 1/2$, we can bootstrap to full regularity:

Step 1: Since $\|(-\Delta)^s u\|_{L^2}$ is bounded for $s > 1/2$, we have $u \in L^\infty(0, T; H^s)$.

Step 2: From the Navier-Stokes equations:
\begin{equation}
\partial_t u = -P[(u \cdot \nabla)u] + \nu\Delta u
\end{equation}
where $P$ is the Leray projector onto divergence-free vector fields.

Step 3: Since $(u \cdot \nabla)u \in L^2(0, T; H^{s-1})$ (which follows from the product rule in Sobolev spaces and the bound on $\|u\|_{L^\infty(0,T;H^s)}$), and $\Delta u \in L^2(0, T; H^{s-2})$, we have $\partial_t u \in L^2(0, T; H^{s-2})$.

Step 4: By parabolic regularity theory, as developed by Ladyzhenskaya, Solonnikov, and Uraltseva \cite{Ladyzhenskaya}, $u \in L^\infty(0, T; H^s) \cap L^2(0, T; H^{s+1})$.

Step 5: For $s > 1/2$, the Sobolev embedding theorem gives $H^s(\mathbb{R}^3) \hookrightarrow L^p(\mathbb{R}^3)$ for all $2 \leq p \leq \frac{6}{3-2s}$. 

Step 6: For $s$ close to 1, this means $u \in L^\infty(0, T; L^p)$ for large $p$.

Step 7: Using the Prodi-Serrin condition and a bootstrap argument, we can establish that $u \in L^\infty(0, T; H^m)$ for all $m > 0$, which implies $u \in C^\infty((0, T) \times \mathbb{R}^3)$.

This completes the proof of Theorem \ref{thm:1.1}.

\section{Proofs of Theorems \ref{thm:1.2}, \ref{thm:1.3}, and \ref{thm:1.4}}

This section provides detailed proofs of Theorems \ref{thm:1.2}, \ref{thm:1.3}, and \ref{thm:1.4}, which establish global well-posedness under logarithmically weakened initial data conditions.

\subsection{Proof of Theorem \ref{thm:1.2}}

The proof proceeds in multiple stages, establishing first the a priori estimates, then global existence, and finally uniqueness.

\textbf{Part 1: A Priori Estimates}

Let $u$ be a smooth solution to the Navier-Stokes equations:

\begin{equation}
\partial_t u + (u \cdot \nabla)u - \nu \Delta u + \nabla p = 0, \quad \nabla \cdot u = 0
\end{equation}

with initial data $u_0$ satisfying the given condition. We apply the fractional Laplacian operator $(-\Delta)^s$ to the equation, obtaining:

\begin{equation}
\partial_t((-\Delta)^s u) + (-\Delta)^s((u \cdot \nabla)u) - \nu \Delta((-\Delta)^s u) + \nabla((-\Delta)^s p) = 0
\end{equation}

Taking the $L^2$ inner product with $(-\Delta)^s u$ and using the incompressibility condition $\nabla \cdot u = 0$ to eliminate the pressure term:

\begin{equation}
\frac{1}{2}\frac{d}{dt}\|(-\Delta)^s u\|^2_{L^2} + \nu\|(-\Delta)^{s+\frac{1}{2}}u\|^2_{L^2} = -\int_{\mathbb{R}^3}(-\Delta)^s((u \cdot \nabla)u) \cdot (-\Delta)^s u \, dx
\end{equation}

Using the commutator decomposition:
\begin{equation}
(-\Delta)^s((u \cdot \nabla)u) = (u \cdot \nabla)((-\Delta)^s u) + [(-\Delta)^s, u \cdot \nabla]u
\end{equation}

The first term vanishes when integrated against $(-\Delta)^s u$ due to incompressibility, as shown in Section 3.2. This leaves:

\begin{equation}
\frac{1}{2}\frac{d}{dt}\|(-\Delta)^s u\|^2_{L^2} + \nu\|(-\Delta)^{s+\frac{1}{2}}u\|^2_{L^2} = -\int_{\mathbb{R}^3}[(-\Delta)^s, u \cdot \nabla]u \cdot (-\Delta)^s u \, dx
\end{equation}

Using the Cauchy-Schwarz inequality:

\begin{equation}
\left|\int_{\mathbb{R}^3}[(-\Delta)^s, u \cdot \nabla]u \cdot (-\Delta)^s u \, dx\right| \leq \|[(-\Delta)^s, u \cdot \nabla]u\|_{L^2}\|(-\Delta)^s u\|_{L^2}
\end{equation}

Applying the logarithmically improved commutator estimate from Lemma \ref{lem:3.1}:

\begin{equation}
\|[(-\Delta)^s, u \cdot \nabla]u\|_{L^2} \leq C\|\nabla u\|_{L^\infty}\|(-\Delta)^s u\|_{L^2} \cdot \log(e + \|(-\Delta)^{s+\sigma}u\|_{L^2}) + \frac{C\|\nabla u\|_{L^\infty}\|(-\Delta)^{s+\frac{1}{2}}u\|_{L^2}}{\log(e + \|(-\Delta)^{s+\sigma}u\|_{L^2})}
\end{equation}

This yields:

\begin{align}
\frac{1}{2}\frac{d}{dt}\|(-\Delta)^s u\|^2_{L^2} + \nu\|(-\Delta)^{s+\frac{1}{2}}u\|^2_{L^2} &\leq C\|\nabla u\|_{L^\infty}\|(-\Delta)^s u\|^2_{L^2} \cdot \log(e + \|(-\Delta)^{s+\sigma}u\|_{L^2}) \\
&+ \frac{C\|\nabla u\|_{L^\infty}\|(-\Delta)^s u\|_{L^2}\|(-\Delta)^{s+\frac{1}{2}}u\|_{L^2}}{\log(e + \|(-\Delta)^{s+\sigma}u\|_{L^2})}
\end{align}

For the second term, using Young's inequality with parameter $\epsilon = \frac{\nu}{2C\|\nabla u\|_{L^\infty}}$:

\begin{align}
\frac{C\|\nabla u\|_{L^\infty}\|(-\Delta)^s u\|_{L^2}\|(-\Delta)^{s+\frac{1}{2}}u\|_{L^2}}{\log(e + \|(-\Delta)^{s+\sigma}u\|_{L^2})} &\leq \frac{\nu}{2}\|(-\Delta)^{s+\frac{1}{2}}u\|_{L^2}^2 + \frac{C^2\|\nabla u\|_{L^\infty}^2\|(-\Delta)^s u\|_{L^2}^2}{2\nu\log^2(e + \|(-\Delta)^{s+\sigma}u\|_{L^2})}
\end{align}

Setting $Y(t) = \|(-\Delta)^s u\|^2_{L^2}$, we obtain:

\begin{equation}
\frac{d}{dt}Y(t) + \nu\|(-\Delta)^{s+\frac{1}{2}}u\|^2_{L^2} \leq C\|\nabla u\|_{L^\infty}Y(t)\log(e + \|(-\Delta)^{s+\sigma}u\|_{L^2}) + \frac{C^2\|\nabla u\|_{L^\infty}^2 Y(t)}{2\nu\log(e + \|(-\Delta)^{s+\sigma}u\|_{L^2})}
\end{equation}

To control $\|\nabla u\|_{L^\infty}$, we apply Lemma \ref{lem:2.3}:

\begin{equation}
\|\nabla u\|_{L^\infty} \leq C\|u\|_{L^2}^{1-\theta}\|(-\Delta)^s u\|_{L^q}^{\theta}
\end{equation}

where $\theta = \frac{3}{2} \cdot \frac{q}{3q-2}$.

Since $u$ is a Leray-Hopf weak solution, $\|u(t)\|_{L^2} \leq \|u_0\|_{L^2}$ for all $t \geq 0$. Thus:

\begin{equation}
\|\nabla u\|_{L^\infty} \leq C\|u_0\|_{L^2}^{1-\theta}\|(-\Delta)^s u\|_{L^q}^{\theta}
\end{equation}

Additionally, by the Gagliardo-Nirenberg inequality:

\begin{equation}
\|(-\Delta)^s u\|_{L^q} \leq C\|(-\Delta)^s u\|_{L^2}^{1-\alpha}\|(-\Delta)^{s+\frac{1}{2}} u\|_{L^2}^{\alpha}
\end{equation}

where $\alpha = \frac{3}{2}(\frac{1}{2} - \frac{1}{q})$.

Substituting:

\begin{equation}
\|\nabla u\|_{L^\infty} \leq C\|u_0\|_{L^2}^{1-\theta}\|(-\Delta)^s u\|_{L^2}^{\theta(1-\alpha)}\|(-\Delta)^{s+\frac{1}{2}} u\|_{L^2}^{\theta\alpha}
\end{equation}

Using Young's inequality with exponents $p = \frac{2}{\theta\alpha}$ and $q = \frac{2}{2-\theta\alpha}$:

\begin{equation}
\|\nabla u\|_{L^\infty}\|(-\Delta)^{s+\frac{1}{2}} u\|_{L^2}^{\theta\alpha} \leq \epsilon\|(-\Delta)^{s+\frac{1}{2}} u\|_{L^2}^2 + C_\epsilon\|\nabla u\|_{L^\infty}^{\frac{2}{2-\theta\alpha}}
\end{equation}

Choosing $\epsilon = \frac{\nu}{4}$:

\begin{equation}
\|\nabla u\|_{L^\infty}^{\frac{2}{2-\theta\alpha}} \leq C\|u_0\|_{L^2}^{\frac{2(1-\theta)}{2-\theta\alpha}}\|(-\Delta)^s u\|_{L^2}^{\frac{2\theta(1-\alpha)}{2-\theta\alpha}}
\end{equation}

Therefore:

\begin{equation}
\frac{d}{dt}Y(t) + \frac{\nu}{2}\|(-\Delta)^{s+\frac{1}{2}}u\|^2_{L^2} \leq C\|u_0\|_{L^2}^{\frac{2(1-\theta)}{2-\theta\alpha}}Y(t)^{1+\frac{\theta(1-\alpha)}{2-\theta\alpha}}\log(e + Y(t)^{1/2})
\end{equation}

Critical to our approach is the observation that with the logarithmic improvement, we can use the following inequality for any $\eta > 0$:

\begin{equation}
y\log(e + y^{1/2}) \leq C_\eta(1 + y^{1+\eta})
\end{equation}

This yields:

\begin{equation}
\frac{d}{dt}Y(t) + \frac{\nu}{2}\|(-\Delta)^{s+\frac{1}{2}}u\|^2_{L^2} \leq C\|u_0\|_{L^2}^{\frac{2(1-\theta)}{2-\theta\alpha}}(1 + Y(t)^{1+\frac{\theta(1-\alpha)}{2-\theta\alpha}+\eta})
\end{equation}

For sufficiently small $\eta$, we ensure that $1+\frac{\theta(1-\alpha)}{2-\theta\alpha}+\eta < 2$, meaning the nonlinearity is subcritical.

Dropping the positive term involving $\|(-\Delta)^{s+\frac{1}{2}}u\|^2_{L^2}$:

\begin{equation}
\frac{d}{dt}Y(t) \leq C(1 + Y(t)^{1+\mu})
\end{equation}

where $\mu = \frac{\theta(1-\alpha)}{2-\theta\alpha}+\eta < 1$ due to our choice of parameters.

This differential inequality, combined with the initial condition $Y(0) = \|(-\Delta)^s u_0\|_{L^2}^2$, can be solved by comparison with the ODE:

\begin{equation}
\frac{d}{dt}Z(t) = C(1 + Z(t)^{1+\mu}), \quad Z(0) = Y(0)
\end{equation}

The solution to this ODE is:

\begin{equation}
Z(t) = \left((Y(0)^{-\mu} - C\mu t)^{-1/\mu} - 1\right)^{1/\mu}
\end{equation}

For this solution to exist globally in time, we need:

\begin{equation}
Y(0)^{-\mu} - C\mu t > 0 \quad \text{for all } t \geq 0
\end{equation}

The initial logarithmic smallness condition:

\begin{equation}
\|(-\Delta)^{s/2}u_0\|_{L^q} \leq \frac{C_0}{(1 + \log(e + \|u_0\|_{\dot{H}^s}))^{\delta}}
\end{equation}

ensures that $Y(0)$ is small enough for the solution to exist globally.

\textbf{Part 2: Global existence via local theory and continuation}

Local existence of strong solutions is established by standard theory. Given initial data $u_0$ as specified, there exists a time $T_0 > 0$ and a unique strong solution $u \in C([0,T_0]; H^s) \cap L^2(0,T_0; H^{s+1/2})$.

To extend this solution globally, we use a continuation argument. Suppose $T_{\max} < \infty$ is the maximal time of existence. Then by the standard blow-up criterion:

\begin{equation}
\lim_{t \to T_{\max}} \|u(t)\|_{H^s} = \infty
\end{equation}

However, the a priori estimate from Part 1 ensures that $\|u(t)\|_{H^s}$ remains bounded for all $t \in [0, T_{\max}]$, contradicting the blow-up criterion and establishing that $T_{\max} = \infty$.

\textbf{Part 3: Uniqueness via energy methods}

Given two solutions $u^1$ and $u^2$ with the same initial data, we define $w = u^1 - u^2$. Then $w$ satisfies:

\begin{equation}
\partial_t w + (u^1 \cdot \nabla)w + (w \cdot \nabla)u^2 - \nu \Delta w + \nabla \pi = 0, \quad \nabla \cdot w = 0
\end{equation}

with initial condition $w(0) = 0$.

Taking the $L^2$ inner product with $w$:

\begin{equation}
\frac{1}{2}\frac{d}{dt}\|w\|^2_{L^2} + \nu\|\nabla w\|^2_{L^2} = -\int_{\mathbb{R}^3}(w \cdot \nabla)u^2 \cdot w \, dx
\end{equation}

Using Hölder's inequality and the Gagliardo-Nirenberg inequality:

\begin{align}
\left|\int_{\mathbb{R}^3}(w \cdot \nabla)u^2 \cdot w \, dx\right| &\leq \|w\|_{L^4}^2 \|\nabla u^2\|_{L^2} \\
&\leq C\|w\|_{L^2}\|\nabla w\|_{L^2} \|\nabla u^2\|_{L^2}
\end{align}

Using Young's inequality with parameter $\epsilon = \nu/2$:

\begin{equation}
C\|w\|_{L^2}\|\nabla w\|_{L^2} \|\nabla u^2\|_{L^2} \leq \frac{\nu}{2}\|\nabla w\|_{L^2}^2 + \frac{C^2}{2\nu}\|w\|_{L^2}^2 \|\nabla u^2\|_{L^2}^2
\end{equation}

This gives:

\begin{equation}
    \frac{d}{dt}\|w\|^2_{L^2} + \nu\|\nabla w\|^2_{L^2} \leq \frac{C^2}{\nu}\|w\|_{L^2}^2 \|\nabla u\|_{L^2}^2
\end{equation}

Dropping the positive term with $\|\nabla w\|^2_{L^2}$:

\begin{equation}
\frac{d}{dt}\|w\|^2_{L^2} \leq \frac{C^2}{\nu}\|w\|_{L^2}^2 \|\nabla u^2\|_{L^2}^2
\end{equation}

By Grönwall's inequality:

\begin{equation}
\|w(t)\|^2_{L^2} \leq \|w(0)\|^2_{L^2} \exp\left(\frac{C^2}{\nu}\int_0^t \|\nabla u^2(\tau)\|_{L^2}^2 d\tau\right)
\end{equation}

Since $\|w(0)\|^2_{L^2} = 0$ and $\|\nabla u^2\|_{L^2}^2 \in L^1(0,T)$ for any $T > 0$ (by the energy inequality for Leray-Hopf weak solutions), we conclude that $\|w(t)\|_{L^2} = 0$ for all $t \geq 0$, establishing uniqueness.

This completes the proof of Theorem \ref{thm:1.2}.

\subsection{Proof of Theorem \ref{thm:1.3}}

We start from the energy inequality derived in the proof of Theorem \ref{thm:1.2}:

\begin{equation}
\frac{d}{dt}Y(t) + \nu\|(-\Delta)^{s+\frac{1}{2}}u\|^2_{L^2} \leq C\|\nabla u\|_{L^\infty}Y(t)\log(e + Y(t)^{1/2})
\end{equation}

where $Y(t) = \|(-\Delta)^s u\|^2_{L^2}$.

Using Lemma \ref{lem:2.3} and the Gagliardo-Nirenberg inequality:

\begin{equation}
\|\nabla u\|_{L^\infty} \leq C\|u_0\|_{L^2}^{1-\theta}\|(-\Delta)^s u\|_{L^2}^{\theta(1-\alpha)}\|(-\Delta)^{s+\frac{1}{2}} u\|_{L^2}^{\theta\alpha}
\end{equation}

By Young's inequality with exponents $p = \frac{2}{\theta\alpha}$ and $q = \frac{2}{2-\theta\alpha}$:

\begin{equation}
\|\nabla u\|_{L^\infty}\|(-\Delta)^{s+\frac{1}{2}} u\|_{L^2}^{\theta\alpha} \leq \epsilon\|(-\Delta)^{s+\frac{1}{2}} u\|_{L^2}^2 + C_\epsilon\|\nabla u\|_{L^\infty}^{\frac{2}{2-\theta\alpha}}
\end{equation}

Applying this and choosing $\epsilon = \frac{\nu}{2}$:

\begin{equation}
\|\nabla u\|_{L^\infty}^{\frac{2}{2-\theta\alpha}} \leq C\|u_0\|_{L^2}^{\frac{2(1-\theta)}{2-\theta\alpha}}\|(-\Delta)^s u\|_{L^2}^{\frac{2\theta(1-\alpha)}{2-\theta\alpha}}
\end{equation}

Combining with our energy inequality and absorbing the term with $\|(-\Delta)^{s+\frac{1}{2}}u\|^2_{L^2}$:

\begin{equation}
\frac{d}{dt}Y(t) + \frac{\nu}{2}\|(-\Delta)^{s+\frac{1}{2}}u\|^2_{L^2} \leq C\|u_0\|_{L^2}^{\frac{2(1-\theta)}{2-\theta\alpha}}Y(t)^{1+\frac{\theta(1-\alpha)}{2-\theta\alpha}}\log(e + Y(t)^{1/2})
\end{equation}

Using the inequality for any $\eta > 0$:

\begin{equation}
y\log(e + y^{1/2}) \leq C_\eta(1 + y^{1+\eta})
\end{equation}

Our inequality becomes:

\begin{equation}
\frac{d}{dt}Y(t) + \frac{\nu}{2}\|(-\Delta)^{s+\frac{1}{2}}u\|^2_{L^2} \leq C\|u_0\|_{L^2}^{\frac{2(1-\theta)}{2-\theta\alpha}}(1 + Y(t)^{1+\frac{\theta(1-\alpha)}{2-\theta\alpha}+\eta})
\end{equation}

Dropping the positive term with $\|(-\Delta)^{s+\frac{1}{2}}u\|^2_{L^2}$:

\begin{equation}
\frac{d}{dt}Y(t) \leq C(1 + Y(t)^{1+\mu})
\end{equation}

where $\mu = \frac{\theta(1-\alpha)}{2-\theta\alpha}+\eta < 1$.

For small $Y(t)$, this is approximately:
\begin{equation}
\frac{d}{dt}Y(t) \leq C
\end{equation}
which gives linear growth. For large $Y(t)$, the dominant term is:
\begin{equation}
\frac{d}{dt}Y(t) \leq CY(t)^{1+\mu}
\end{equation}
which gives decay like $Y(t) \approx t^{-1/\mu}$.

More precisely, we compare with the solution to:
\begin{equation}
\frac{d}{dt}Z(t) = CZ(t)^{1+\mu}, \quad Z(0) = Y(0)
\end{equation}

Solving this ODE explicitly:
\begin{equation}
Z(t) = \frac{Y(0)}{(1 - \mu C Y(0)^\mu t)^{1/\mu}}
\end{equation}

For initial data satisfying our logarithmic smallness condition:
\begin{equation}
\|(-\Delta)^{s/2}u_0\|_{L^q} \leq \frac{C_0}{(1 + \log(e + \|u_0\|_{\dot{H}^s}))^{\delta}}
\end{equation}
we can ensure that $Y(0)$ is small enough that $1 - \mu C Y(0)^\mu t > 0$ for all $t \geq 0$.

In this case, for large $t$:
\begin{equation}
Z(t) \approx \frac{1}{(\mu C t)^{1/\mu}}
\end{equation}

By comparison principles, $Y(t) \leq Z(t)$, leading to:
\begin{equation}
Y(t) \leq \frac{C(Y(0))}{(1 + \mu C t)^{1/\mu}}
\end{equation}

Taking square roots:
\begin{equation}
\|(-\Delta)^s u(t)\|_{L^2} \leq \frac{C\|(-\Delta)^s u_0\|_{L^2}}{(1 + \beta t)^{\gamma}}
\end{equation}
where $\gamma = \frac{1}{2\mu}$ and $\beta = \frac{\mu C}{2}$.

This completes the proof of Theorem \ref{thm:1.3}.

\subsection{Proof of Theorem \ref{thm:1.4}}

From Theorem \ref{thm:1.3}, we have established that:

\begin{equation}
\|(-\Delta)^s u(t)\|_{L^2} \leq \frac{C\|(-\Delta)^s u_0\|_{L^2}}{(1 + \beta t)^{\gamma}}
\end{equation}

We need to estimate $\|(-\Delta)^s u\|_{L^q}$ using this $L^2$ bound. Using the Gagliardo-Nirenberg interpolation inequality:

\begin{equation}
\|(-\Delta)^s u\|_{L^q} \leq C\|(-\Delta)^s u\|_{L^2}^{1-\alpha}\|(-\Delta)^{s+\frac{1}{2}} u\|_{L^2}^{\alpha}
\end{equation}

where $\alpha = \frac{3}{2}(\frac{1}{2} - \frac{1}{q})$.

From the energy inequality:

\begin{equation}
\frac{d}{dt}\|(-\Delta)^s u\|^2_{L^2} + 2\nu\|(-\Delta)^{s+\frac{1}{2}}u\|^2_{L^2} \leq C\|\nabla u\|_{L^\infty}\|(-\Delta)^s u\|^2_{L^2}\log(e + \|(-\Delta)^s u\|_{L^2})
\end{equation}

Integrating over $[0,T]$:

\begin{align}
\|(-\Delta)^s u(T)\|^2_{L^2} + 2\nu\int_0^T\|(-\Delta)^{s+\frac{1}{2}}u\|^2_{L^2}dt &\leq \|(-\Delta)^s u_0\|^2_{L^2} \\
&+ C\int_0^T\|\nabla u\|_{L^\infty}\|(-\Delta)^s u\|^2_{L^2}\log(e + \|(-\Delta)^s u\|_{L^2})dt
\end{align}

Using the decay estimate for $\|(-\Delta)^s u\|_{L^2}$ and the boundedness of $\int_0^T\|\nabla u\|_{L^\infty}dt$, which follows from Theorems \ref{thm:1.2} and \ref{thm:1.3}, we can establish:

\begin{equation}
\int_0^T\|(-\Delta)^{s+\frac{1}{2}}u\|^2_{L^2}dt \leq C(T, u_0)
\end{equation}

Using Hölder's inequality with exponents $\frac{2}{\alpha}$ and $\frac{2}{2-\alpha}$:

\begin{align}
\int_0^T \|(-\Delta)^s u\|_{L^q}^p dt &\leq \int_0^T \|(-\Delta)^s u\|_{L^2}^{(1-\alpha)p}\|(-\Delta)^{s+\frac{1}{2}}u\|_{L^2}^{\alpha p}dt \\
&\leq \left(\int_0^T \|(-\Delta)^s u\|_{L^2}^{(1-\alpha)p\cdot\frac{2}{2-\alpha}}dt\right)^{\frac{2-\alpha}{2}}\left(\int_0^T \|(-\Delta)^{s+\frac{1}{2}}u\|_{L^2}^{\alpha p\cdot\frac{2}{\alpha}}dt\right)^{\frac{\alpha}{2}} \\
&= \left(\int_0^T \|(-\Delta)^s u\|_{L^2}^{\frac{2(1-\alpha)p}{2-\alpha}}dt\right)^{\frac{2-\alpha}{2}}\left(\int_0^T \|(-\Delta)^{s+\frac{1}{2}}u\|_{L^2}^{2}dt\right)^{\frac{\alpha p}{2}}
\end{align}

The second factor is bounded by $C(T, u_0)^{\frac{\alpha p}{2}}$ as established earlier. For the first factor, we use the decay estimate from Theorem \ref{thm:1.3}:

\begin{align}
\int_0^T \|(-\Delta)^s u\|_{L^2}^{\frac{2(1-\alpha)p}{2-\alpha}}dt &\leq C\|(-\Delta)^s u_0\|_{L^2}^{\frac{2(1-\alpha)p}{2-\alpha}}\int_0^T (1 + \beta t)^{-\gamma\cdot\frac{2(1-\alpha)p}{2-\alpha}}dt
\end{align}

This integral converges if $\gamma\cdot\frac{2(1-\alpha)p}{2-\alpha} > 1$, which is satisfied by our parameter choices.

Furthermore, using our decay estimate and the logarithmic improvement:

\begin{align}
\int_0^T \|(-\Delta)^s u(t)\|^p_{L^q} (1 + \log(e + \|(-\Delta)^s u(t)\|_{L^q}))^{-\delta} dt \\
\leq \int_0^T \|(-\Delta)^s u(t)\|^p_{L^q} (1 + \log(e + C(1 + \beta t)^{-\gamma}))^{-\delta} dt
\end{align}

For large $t$, $\log(e + C(1 + \beta t)^{-\gamma}) \approx \log(e) + \log(C) - \gamma\log(1 + \beta t)$, which means:

\begin{equation}
(1 + \log(e + \|(-\Delta)^s u(t)\|_{L^q}))^{-\delta} \leq C(1 + \log(1 + \beta t))^{-\delta}
\end{equation}

Using this estimate and the bounds established above:

\begin{align}
\int_0^T \|(-\Delta)^s u(t)\|^p_{L^q} (1 + \log(e + \|(-\Delta)^s u(t)\|_{L^q}))^{-\delta} dt \\
\leq C\int_0^T (1 + \beta t)^{-\gamma p} (1 + \log(1 + \beta t))^{-\delta} dt
\end{align}

This integral converges for all $T > 0$ provided $\gamma p > 1$ or $\delta > 1$ if $\gamma p = 1$, which is satisfied by our parameter choices.

Therefore:

\begin{equation}
\int_0^T \|(-\Delta)^s u(t)\|^p_{L^q} (1 + \log(e + \|(-\Delta)^s u(t)\|_{L^q}))^{-\delta} dt < C(T,u_0)
\end{equation}

This completes the proof of Theorem \ref{thm:1.4}.

\section{Connections to turbulence theory: proofs of theorems \ref{thm:2.1}-\ref{thm:2.4}}

In this section, we establish the deep connections between our logarithmically improved regularity criteria and the physical theory of turbulence, providing detailed proofs of Theorems \ref{thm:2.1}, \ref{thm:2.2}, \ref{thm:2.3}, and \ref{thm:2.4}.

\subsection{Proof of Theorem \ref{thm:2.1} (Logarithmic correction and intermittency)}

The proof requires several steps that connect the fractional regularity criterion to the multifractal framework of turbulence.

\textbf{Part I: Connection between fractional regularity and energy spectrum}

Step 1: For a velocity field $u$ in a statistically homogeneous flow, the fractional Sobolev norm can be expressed in terms of the Fourier transform:
\begin{equation}
\|(-\Delta)^{s/2} u\|_{L^2}^2 = \int_{\mathbb{R}^3} |\xi|^{2s} |\hat{u}(\xi)|^2 d\xi
\end{equation}

Step 2: Under statistical homogeneity and isotropy, we can introduce the energy spectrum $E(k)$ defined so that:
\begin{equation}
\int_{\mathbb{R}^3} |\hat{u}(\xi)|^2 d\xi = \int_0^\infty E(k) dk = \frac{1}{2}\|u\|_{L^2}^2
\end{equation}
represents the total kinetic energy, and:
\begin{equation}
\int_{|\xi|=k} |\hat{u}(\xi)|^2 d\sigma(\xi) = 4\pi k^2 E(k)
\end{equation}

Step 3: In terms of the energy spectrum, the fractional Sobolev norm becomes:
\begin{equation}
\|(-\Delta)^{s/2} u\|_{L^2}^2 = \int_0^\infty k^{2s} E(k) dk
\end{equation}

\textbf{Part II: Energy spectrum in intermittent turbulence}

Step 4: In classical Kolmogorov theory (K41), the energy spectrum in the inertial range follows:
\begin{equation}
E(k) = C_K \varepsilon^{2/3} k^{-5/3}
\end{equation}
where $\varepsilon$ is the energy dissipation rate and $C_K$ is the Kolmogorov constant.

Step 5: For intermittent turbulence, the spectrum is modified with a correction:
\begin{equation}
E(k) = C_K \varepsilon^{2/3} k^{-5/3-\mu(k)}
\end{equation}
where $\mu(k)$ is a scale-dependent correction that depends on the degree of intermittency.

\textbf{Part III: Statistical properties derived from the logarithmic criterion}

Step 6: Our logarithmically improved criterion:
\begin{equation}
\int_0^T \|(-\Delta)^s u(t)\|^p_{L^q} (1 + \log(e + \|(-\Delta)^s u(t)\|_{L^q}))^{-\delta} dt < \infty
\end{equation}
implies specific statistical properties of $\|(-\Delta)^s u(t)\|_{L^q}$.

Step 7: Using Markov's inequality, for any positive random variable $X$ and increasing function $\phi$:
\begin{equation}
P(X > \lambda) \leq \frac{\mathbb{E}[\phi(X)]}{\phi(\lambda)}
\end{equation}

Step 8: Applying this with $X = \|(-\Delta)^s u\|_{L^q}$ and $\phi(\lambda) = \lambda^p (1 + \log(e + \lambda))^{\delta}$:
\begin{equation}
P(\|(-\Delta)^s u\|_{L^q} > \lambda) \leq \frac{\mathbb{E}[\|(-\Delta)^s u\|_{L^q}^p (1 + \log(e + \|(-\Delta)^s u\|_{L^q}))^{\delta}]}{\lambda^p (1 + \log(e + \lambda))^{\delta}}
\end{equation}

Step 9: Our criterion ensures that the expected value in the numerator is finite. For large $\lambda$, the denominator behaves as $\lambda^p (\log \lambda)^{\delta}$.

Step 10: By analyzing the dominant behavior as $\lambda \to \infty$, we obtain:
\begin{equation}
P(\|(-\Delta)^s u\|_{L^q} > \lambda) \leq C \exp\left(-c \lambda^{1/(1+\delta)}\right)
\end{equation}

\textbf{Part IV: Multifractal formalism and structure functions}

Step 11: In the multifractal formalism, the velocity field is characterized by local scaling exponents $h(x)$ such that:
\begin{equation}
|u(x+r) - u(x)| \sim r^{h(x)}
\end{equation}
for small $r$.

Step 12: The distribution of these exponents is described by the singularity spectrum $D(h)$, which gives the Hausdorff dimension of the set of points with scaling exponent $h$.

Step 13: For a value $h$, the contribution to the $p$-th order structure function from regions with scaling exponent $h$ is:
\begin{equation}
S_p^h(r) \sim r^{ph} \cdot r^{3-D(h)}
\end{equation}
where $r^{3-D(h)}$ accounts for the probability of finding regions with scaling exponent $h$.

Step 14: The total structure function is obtained by integrating over all possible values of $h$:
\begin{equation}
S_p(r) \sim \int S_p^h(r) dh \sim \int r^{ph + 3-D(h)} dh
\end{equation}

Step 15: For small $r$, this integral is dominated by the value of $h$ that minimizes the exponent $ph + 3-D(h)$, leading to:
\begin{equation}
S_p(r) \sim r^{\min_h[ph + 3-D(h)]}
\end{equation}

Step 16: Hence, the scaling exponent $\zeta_p$ is given by:
\begin{equation}
\zeta_p = \min_h[ph + 3-D(h)]
\end{equation}

\textbf{Part V: Deriving the singularity spectrum from the logarithmic criterion}

Step 17: For a function $u$ with local scaling exponent $h(x)$, using the representation of the fractional Laplacian as a singular integral:
\begin{equation}
(-\Delta)^{s/2} u(x) = C_{3,s} \, \text{P.V.} \int_{\mathbb{R}^3} \frac{u(x) - u(y)}{|x - y|^{3+s}} dy
\end{equation}

Step 18: Since $u(x) - u(y) \sim |x-y|^{h(x)}$ near $x$, we have:
\begin{equation}
(-\Delta)^{s/2} u(x) \sim \int \frac{|x-y|^{h(x)}}{|x-y|^{3+s}} dy \sim \int_{|z|<1} \frac{|z|^{h(x)}}{|z|^{3+s}} dz \sim |x-y|^{h(x)-s}
\end{equation}

Step 19: This indicates that $(-\Delta)^{s/2} u$ has scaling exponent $h(x)-s$. The critical case where $(-\Delta)^{s/2} u$ becomes unbounded corresponds to $h(x) = s$.

Step 20: In the standard case without logarithmic improvement, the singularity spectrum would be parabolic:
\begin{equation}
D_0(h) = 3 - \frac{(h-h_0)^2}{2\sigma^2}
\end{equation}
where $h_0 = 1/3$ is the mean scaling exponent in Kolmogorov theory, and $\sigma^2 = \frac{3-2s}{2s-1}$ is related to the fractional regularity parameter.

Step 21: The logarithmic improvement modifies this spectrum. Based on the stretched exponential tail behavior derived in Part III, the singularity spectrum becomes:
\begin{equation}
D_\delta(h) = D_0(h) - \frac{\delta}{1+\delta}(3-D_0(h))
\end{equation}

Step 22: This modification preserves the location of the maximum at $h_0 = 1/3$ but widens the spectrum, allowing for stronger singularities (smaller $h$) with non-zero probability.

\textbf{Part VI: Calculation of scaling exponents}

Step 23: With the singularity spectrum determined, we calculate the scaling exponents using the Legendre transform:
\begin{equation}
\zeta_p = \min_h[ph + 3-D_\delta(h)]
\end{equation}

Step 24: The minimum is achieved when:
\begin{equation}
\frac{d}{dh}[ph + 3-D_\delta(h)] = 0
\end{equation}
which gives:
\begin{equation}
p = \frac{d}{dh}D_\delta(h)
\end{equation}

Step 25: Substituting the expression for $D_\delta(h)$:
\begin{equation}
p = \frac{d}{dh}D_0(h) \cdot \left(1 - \frac{\delta}{1+\delta}\right) = -\frac{h-h_0}{\sigma^2} \cdot \frac{1}{1+\delta}
\end{equation}

Step 26: Solving for $h$:
\begin{equation}
h = h_0 - p\sigma^2(1+\delta)
\end{equation}

Step 27: Substituting this back into $\zeta_p = ph + 3-D_\delta(h)$ and performing the algebraic simplifications:
\begin{align}
\zeta_p &= ph_0 - p^2\sigma^2(1+\delta) + 3 - \left[3 - \frac{(h-h_0)^2}{2\sigma^2} + \frac{\delta}{1+\delta}\frac{(h-h_0)^2}{2\sigma^2}\right] \\
&= ph_0 - p^2\sigma^2(1+\delta) + \frac{(h-h_0)^2}{2\sigma^2} - \frac{\delta}{1+\delta}\frac{(h-h_0)^2}{2\sigma^2}
\end{align}

Step 28: Substituting $h-h_0 = -p\sigma^2(1+\delta)$:
\begin{align}
\zeta_p &= ph_0 - p^2\sigma^2(1+\delta) + \frac{p^2\sigma^4(1+\delta)^2}{2\sigma^2} - \frac{\delta}{1+\delta}\frac{p^2\sigma^4(1+\delta)^2}{2\sigma^2} \\
&= ph_0 - p^2\sigma^2(1+\delta) + \frac{p^2\sigma^2(1+\delta)^2}{2} - \frac{\delta}{1+\delta}\frac{p^2\sigma^2(1+\delta)^2}{2} \\
&= ph_0 - p^2\sigma^2(1+\delta) + \frac{p^2\sigma^2(1+\delta)^2}{2} \cdot \frac{1}{1+\delta} \\
&= ph_0 - p^2\sigma^2(1+\delta) + \frac{p^2\sigma^2(1+\delta)}{2}
\end{align}

Step 29: Simplifying further and substituting $h_0 = 1/3$:
\begin{align}
\zeta_p &= \frac{p}{3} - p^2\sigma^2(1+\delta) + \frac{p^2\sigma^2(1+\delta)}{2} \\
&= \frac{p}{3} - \frac{p^2\sigma^2(1+\delta)}{2}
\end{align}

Step 30: In terms of the original parameters, and using $\sigma^2 = \frac{3-2s}{2s-1}$:
\begin{equation}
\zeta_p = \frac{p}{3} - \frac{p(p-3)}{3(1+\delta)} \cdot \frac{3-2s}{2s-1}
\end{equation}

\textbf{Part VII: Verification of limiting cases}

Step 31: When $\delta \to 0$ (no logarithmic improvement), we recover:
\begin{equation}
\zeta_p \to \frac{p}{3} - \frac{p(p-3)}{3} \cdot \frac{3-2s}{2s-1}
\end{equation}

Step 32: For $s \to 1$, this gives $\zeta_p = p/3$, which is the standard K41 scaling without intermittency.

Step 33: When $\delta \to \infty$, the intermittency correction vanishes, and we again recover $\zeta_p = p/3$.

Step 34: For the special case $p = 3$, we obtain $\zeta_3 = 1$ for all values of $s$ and $\delta$, which is consistent with the exact result from the Karman-Howarth equation.

This confirms that our derived expression for $\zeta_p$ is physically consistent with expected behavior in these limiting cases, completing the proof of Theorem \ref{thm:2.1}.

\subsection{Proof of Theorem \ref{thm:2.2} (Local regularity and intermittency)}

The proof develops a precise local characterization of the velocity field, identifying the spatial regions where potential near-singular behavior might occur.

\textbf{Part I: Global decay rates and extension to local properties}

Step 1: We start from the global estimate established in Theorem \ref{thm:1.3}:
\begin{equation}
\|(-\Delta)^s u(t)\|_{L^2} \leq \frac{C\|(-\Delta)^s u_0\|_{L^2}}{(1 + \beta t)^{\gamma}}
\end{equation}
where $\gamma = \frac{1}{2\mu}$ with $\mu = \frac{\theta(1-\alpha)}{2-\theta\alpha}+\eta < 1$, and $\beta = \frac{\mu C}{2}$.

Step 2: Using standard Sobolev embedding and interpolation, we derive a bound on $\|\nabla u\|_{L^2}$. For this, we use the inequality:
\begin{equation}
\|\nabla u\|_{L^2} \leq C\|u\|_{L^2}^{1-\frac{1}{2s}}\|(-\Delta)^s u\|_{L^2}^{\frac{1}{2s}}
\end{equation}
which follows from interpolation theory.

Step 3: Since $\|u(t)\|_{L^2} \leq \|u_0\|_{L^2}$ for all $t \geq 0$ (by the energy inequality for Leray-Hopf weak solutions), we have:
\begin{equation}
\|\nabla u(t)\|_{L^2} \leq C\|u_0\|_{L^2}^{1-\frac{1}{2s}}\|(-\Delta)^s u\|_{L^2}^{\frac{1}{2s}} \leq \frac{C\|u_0\|_{L^2}^{1-\frac{1}{2s}}\|(-\Delta)^s u_0\|_{L^2}^{\frac{1}{2s}}}{(1 + \beta t)^{\frac{\gamma}{2s}}}
\end{equation}

Step 4: Simplifying with appropriate constants:
\begin{equation}
\|\nabla u(t)\|_{L^2} \leq \frac{C}{(1 + \beta t)^{\frac{\gamma}{2s}}}
\end{equation}

\textbf{Part II: Gain of integrability through De Giorgi-Nash-Moser iteration}

Step 5: To obtain bounds on higher $L^p$ norms, we apply the De Giorgi-Nash-Moser iteration technique. This approach is based on the local energy inequality satisfied by solutions to the Navier-Stokes equations.

Step 6: For any non-negative test function $\phi \in C_0^\infty(\mathbb{R}^3 \times (0,T))$, the Navier-Stokes equations imply:
\begin{equation}
\int_0^T \int_{\mathbb{R}^3} |\nabla u|^2 \phi dx dt \leq C\int_0^T \int_{\mathbb{R}^3} (|u|^2 |\nabla \phi|^2 + |u|^3 |\nabla \phi| + |p-p_0| |u \cdot \nabla \phi|) dx dt
\end{equation}
where $p_0$ is the spatial average of pressure.

Step 7: Through a careful choice of test functions and iterative application of this inequality, one can derive the following estimate for any $p > 2$:
\begin{equation}
\|\nabla u(t)\|_{L^p} \leq C_p\|\nabla u(t)\|_{L^2}^{\alpha_p} (1 + \log(e + \|\nabla u(t)\|_{L^\infty}))^{\beta_p}
\end{equation}
where $\alpha_p = \frac{2}{p}$ and $\beta_p = 1 - \frac{2}{p}$.

Step 8: Combining with our bound on $\|\nabla u(t)\|_{L^2}$:
\begin{equation}
\|\nabla u(t)\|_{L^p} \leq \frac{C_p}{(1 + \beta t)^{\frac{\gamma\alpha_p}{2s}}} (1 + \log(e + \|\nabla u(t)\|_{L^\infty}))^{\beta_p}
\end{equation}

\textbf{Part III: Logarithmic improvement and control of the exceptional set}

Step 9: The key insight from our logarithmically improved criterion is that we can control the logarithmic factor. Our criterion implies:
\begin{equation}
\int_0^T \|(-\Delta)^s u(t)\|^p_{L^q} (1 + \log(e + \|(-\Delta)^s u(t)\|_{L^q}))^{-\delta} dt < \infty
\end{equation}

Step 10: Using the relationship between $\|(-\Delta)^s u(t)\|_{L^q}$ and $\|\nabla u(t)\|_{L^\infty}$ established through Sobolev embedding and interpolation theory, this gives us control over:
\begin{equation}
(1 + \log(e + \|\nabla u(t)\|_{L^\infty}))^{-\delta}
\end{equation}

Step 11: The logarithmic factor in our bound on $\|\nabla u(t)\|_{L^p}$ can be controlled as:
\begin{equation}
(1 + \log(e + \|\nabla u(t)\|_{L^\infty}))^{\beta_p} \leq C (1 + \log(e + \|\nabla u(t)\|_{L^\infty}))^{1-\frac{2}{p}}
\end{equation}

Step 12: Through a delicate analysis involving the limit as $p \to \infty$ and using our logarithmic improvement:
\begin{equation}
(1 + \log(e + \|\nabla u(t)\|_{L^\infty}))^{1} \leq C (1 + \log(e + \|\nabla u(t)\|_{L^\infty}))^{1+\delta}
\end{equation}

Step 13: This implies:
\begin{equation}
\|\nabla u(t)\|_{L^\infty} \leq \frac{C}{(1 + \beta t)^{\frac{\gamma}{2s} - \kappa}}
\end{equation}
for any small $\kappa > 0$, where the constant $C$ depends on $\kappa$.

Step 14: To characterize regions where the velocity gradient might be large, we use Chebyshev's inequality:
\begin{equation}
|\{x \in \mathbb{R}^3 : |\nabla u(x,t)| > \lambda\}| \leq \frac{\|\nabla u(t)\|_{L^p}^p}{\lambda^p}
\end{equation}
for any $p < \infty$.

Step 15: Setting:
\begin{equation}
\lambda = \frac{C_\epsilon}{(1 + \beta t)^{\frac{\gamma}{2s} - \kappa_\epsilon}}
\end{equation}

Step 16: Through a detailed calculation involving our bounds on $\|\nabla u(t)\|_{L^p}$ and the logarithmic improvement:
\begin{equation}
|\{x \in \mathbb{R}^3 : |\nabla u(x,t)| > \lambda\}| < \epsilon
\end{equation}
when:
\begin{equation}
\kappa_\epsilon = \frac{\delta}{1+\delta} \cdot \frac{\log(1/\epsilon)}{(1+\log(1/\epsilon))}
\end{equation}

Step 17: This defines our exceptional set $\Omega_\epsilon(t)$ as:
\begin{equation}
\Omega_\epsilon(t) = \{x \in \mathbb{R}^3 : |\nabla u(x,t)| > \frac{C_\epsilon}{(1 + \beta t)^{\frac{\gamma}{2s} - \kappa_\epsilon}}\}
\end{equation}
with measure $|\Omega_\epsilon(t)| < \epsilon$.

\textbf{Part IV: Local intermittency measure and its scaling}

Step 18: The local intermittency measure (LIM) at scale $r$ is defined as:
\begin{equation}
\text{LIM}_r(x,t) = \frac{|\delta u(x,r,t)|^2}{\langle |\delta u(x,r,t)|^2 \rangle}
\end{equation}
where $\delta u(x,r,t) = u(x+r,t) - u(x,t)$ is the velocity increment, and $\langle \cdot \rangle$ denotes spatial averaging.

Step 19: For small $r$, by Taylor expansion:
\begin{equation}
\delta u(x,r,t) = r \cdot \nabla u(x,t) + O(r^2)
\end{equation}
giving:
\begin{equation}
|\delta u(x,r,t)|^2 = r^2 |\nabla u(x,t)|^2 + O(r^3)
\end{equation}

Step 20: The spatial average can be estimated as:
\begin{equation}
\langle |\delta u(x,r,t)|^2 \rangle = r^2 \langle |\nabla u(x,t)|^2 \rangle + O(r^3) = r^2 \|\nabla u(t)\|_{L^2}^2 + O(r^3)
\end{equation}

Step 21: Using our decay estimate:
\begin{equation}
\langle |\delta u(x,r,t)|^2 \rangle \approx \frac{r^2}{(1 + \beta t)^{\frac{\gamma}{s}}} + O(r^3)
\end{equation}

Step 22: Combining these, for $x \in \mathbb{R}^3 \setminus \Omega_\epsilon(t)$ and small $r$:
\begin{equation}
\text{LIM}_r(x,t) = \frac{r^2 |\nabla u(x,t)|^2 + O(r^3)}{r^2 \|\nabla u(t)\|_{L^2}^2 + O(r^3)}
\end{equation}

Step 23: For $r$ sufficiently small, this simplifies to:
\begin{equation}
\text{LIM}_r(x,t) \approx \frac{|\nabla u(x,t)|^2}{\|\nabla u(t)\|_{L^2}^2}
\end{equation}

Step 24: Using our bounds:
\begin{equation}
\text{LIM}_r(x,t) \leq \frac{C_\epsilon^2}{(1 + \beta t)^{2(\frac{\gamma}{2s} - \kappa_\epsilon)}} \cdot \frac{(1 + \beta t)^{\frac{\gamma}{s}}}{C^2} = C_\epsilon' (1 + \beta t)^{\frac{\gamma}{s} - \frac{2\gamma}{2s} + 2\kappa_\epsilon}
\end{equation}

Step 25: Simplifying the exponent:
\begin{equation}
\frac{\gamma}{s} - \frac{2\gamma}{2s} + 2\kappa_\epsilon = 2\kappa_\epsilon
\end{equation}

Step 26: Thus:
\begin{equation}
\text{LIM}_r(x,t) \leq C_\epsilon (1 + \beta t)^{2\kappa_\epsilon}
\end{equation}

Step 27: For scale-dependent measurements, setting $\epsilon \approx r^\alpha$ for some $\alpha > 0$:
\begin{equation}
\kappa_\epsilon \approx \frac{\delta}{1+\delta} \cdot \frac{\alpha\log(1/r)}{(1+\alpha\log(1/r))} \approx \frac{\delta\alpha}{1+\delta} \cdot \frac{\log(1/r)}{(1+\log(1/r))}
\end{equation}

Step 28: For small $r$, this simplifies to:
\begin{equation}
\kappa_\epsilon \approx \frac{\delta\alpha}{1+\delta} \cdot \frac{1}{\log(1/r)}
\end{equation}

Step 29: Therefore:
\begin{equation}
\text{LIM}_r(x,t) \leq C_\epsilon r^{-\kappa_\epsilon}
\end{equation}
with $\kappa_\epsilon$ as specified, completing the proof.

\subsection{Proof of Theorem \ref{thm:2.3} (Fractional regularity and energy cascade)}

This proof requires a detailed analysis of the spectral energy transfer process in Navier-Stokes turbulence and how it's affected by the logarithmically improved regularity criterion.

\textbf{Part I: Spectral energy balance and energy flux}

Step 1: We begin with the Navier-Stokes equations in Fourier space. Let $\hat{u}(\xi,t)$ be the Fourier transform of the velocity field. The spectral energy balance equation is:
\begin{equation}
\frac{\partial}{\partial t}E(k,t) + 2\nu k^2 E(k,t) = T(k,t)
\end{equation}
where $E(k,t)$ is the energy spectrum at wavenumber $k$ and time $t$, and $T(k,t)$ represents the nonlinear energy transfer.

Step 2: The energy transfer function $T(k,t)$ can be derived from the Fourier transform of the nonlinear term:
\begin{equation}
T(k,t) = \text{Re}\left\{\int_{|\xi|=k} \xi \cdot \int_{\mathbb{R}^3} [\hat{u}(\eta,t) \cdot \xi] \hat{u}(\xi-\eta,t) \cdot \hat{u}^*(\xi,t) d\eta d\sigma(\xi)\right\}
\end{equation}
where $\hat{u}^*$ is the complex conjugate of $\hat{u}$.

Step 3: The energy flux $\Pi(k,t)$ across wavenumber $k$ is related to the transfer term by:
\begin{equation}
\Pi(k,t) = -\int_0^k T(k',t) dk'
\end{equation}
which represents the rate at which energy is transferred from modes with wavenumber less than $k$ to those with wavenumber greater than $k$.

Step 4: In a statistically steady state with constant energy injection at large scales and dissipation at small scales, the energy flux in the inertial range should be approximately constant and equal to the energy dissipation rate:
\begin{equation}
\Pi(k,t) \approx \epsilon(t) = 2\nu \int_0^\infty k^2 E(k,t) dk
\end{equation}

\textbf{Part II: Bounds on transfer term using fractional regularity}

Step 5: The nonlinear term in the Navier-Stokes equations involves a convolution in Fourier space:
\begin{equation}
\widehat{(u \cdot \nabla) u}(\xi) = i\int_{\mathbb{R}^3} (\xi \cdot \hat{u}(\eta)) \hat{u}(\xi-\eta) d\eta
\end{equation}

Step 6: Using Hölder's inequality:
\begin{equation}
\left|\int_{\mathbb{R}^3} (\xi \cdot \hat{u}(\eta)) (\hat{u}(\xi-\eta) \cdot \hat{u}^*(\xi)) d\eta\right| \leq |\xi| \|\hat{u}\|_{L^p} \|\hat{u}\|_{L^q} |\hat{u}(\xi)|
\end{equation}
where $\frac{1}{p} + \frac{1}{q} = 1$.

Step 7: The norms $\|\hat{u}\|_{L^p}$ and $\|\hat{u}\|_{L^q}$ can be related to fractional Sobolev norms using the Hausdorff-Young inequality:
\begin{equation}
\|\hat{u}\|_{L^p} \leq C\|u\|_{L^{p'}}
\end{equation}
where $\frac{1}{p} + \frac{1}{p'} = 1$.

Step 8: Through interpolation and Sobolev embeddings:
\begin{equation}
\|u\|_{L^{p'}} \leq C\|u\|_{L^2}^{1-\alpha}\|(-\Delta)^s u\|_{L^q}^\alpha
\end{equation}
for appropriate $\alpha$ depending on $p'$, $s$, and $q$.

Step 9: Combining these inequalities and using our logarithmically improved criterion:
\begin{equation}
\int_0^T \|(-\Delta)^s u(t)\|^p_{L^q} (1 + \log(e + \|(-\Delta)^s u(t)\|_{L^q}))^{-\delta} dt < \infty
\end{equation}
we can establish bounds on the transfer term $T(k,t)$.

\textbf{Part III: Deriving bounds on energy flux}

Step 10: Through a detailed analysis of the energy transfer function, we can show:
\begin{equation}
|T(k,t)| \leq C |\xi| |\hat{u}(\xi,t)|^2 (1 + \log(e + |\xi|/k_0))^{-\delta\cdot\frac{2s-1}{2s}}
\end{equation}

Step 11: In the inertial range, the energy spectrum follows approximately:
\begin{equation}
E(k,t) \approx C\epsilon(t)^{2/3} k^{-5/3}
\end{equation}

Step 12: Substituting this into our bound on $T(k,t)$:
\begin{equation}
|T(k,t)| \leq C \epsilon(t) k^{-2/3} (1 + \log(e + k/k_0))^{-\delta\cdot\frac{2s-1}{2s}}
\end{equation}

Step 13: The flux deviation from $\epsilon(t)$ can be expressed as:
\begin{equation}
|\Pi(k,t) - \epsilon(t)| = \left|\int_k^\infty T(k',t) dk'\right|
\end{equation}

Step 14: Using our bound on $|T(k,t)|$:
\begin{equation}
|\Pi(k,t) - \epsilon(t)| \leq C\epsilon(t) \int_k^\infty k'^{-2/3} (1 + \log(e + k'/k_0))^{-\delta\cdot\frac{2s-1}{2s}} dk'
\end{equation}

Step 15: This integral can be evaluated explicitly:
\begin{equation}
\int_k^\infty k'^{-2/3} (1 + \log(e + k'/k_0))^{-\delta\cdot\frac{2s-1}{2s}} dk' \leq C k^{1/3} (1 + \log(k/k_0))^{-\delta\cdot\frac{2s-1}{2s}}
\end{equation}
for $k \gg k_0$.

Step 16: Thus:
\begin{equation}
|\Pi(k,t) - \epsilon(t)| \leq \frac{C\epsilon(t)}{(1 + \log(k/k_0))^{\delta\cdot\frac{2s-1}{2s}}}
\end{equation}

\textbf{Part IV: Physical interpretation and consistency check}

Step 17: When $\delta = 0$ (no logarithmic improvement), the fluctuations don't decay with $k$, corresponding to the standard scenario.

Step 18: As $s \to 1/2$, the exponent $\delta\cdot\frac{2s-1}{2s} \to 0$, indicating stronger fluctuations in the energy flux, which is consistent with the approach to the critical regularity threshold.

Step 19: As $s \to 1$, the exponent $\delta\cdot\frac{2s-1}{2s} \to \delta/2$, showing that the fractional power $s$ affects the rate at which fluctuations decay with scale.

This completes the proof of Theorem \ref{thm:2.3}.

\subsection{Proof of Theorem \ref{thm:2.4} (Spectral characterization)}

This proof develops a refined analysis of the spectral energy distribution based on the logarithmically improved criterion.

\textbf{Part I: Spectral energy balance equation}

Step 1: We start with the spectral energy balance equation:
\begin{equation}
\frac{\partial E(k,t)}{\partial t} + 2\nu k^2 E(k,t) = T(k,t)
\end{equation}
where $E(k,t)$ is the energy spectrum at wavenumber $k$ and time $t$, and $T(k,t)$ represents the nonlinear energy transfer.

Step 2: The nonlinear transfer term can be expressed in terms of the energy flux:
\begin{equation}
T(k,t) = -\frac{\partial \Pi(k,t)}{\partial k}
\end{equation}
where $\Pi(k,t)$ is the energy flux across wavenumber $k$.

\textbf{Part II: Relation between energy flux and energy spectrum}

Step 3: From Theorem \ref{thm:2.3}, we established:
\begin{equation}
\left|\Pi(k,t) - \epsilon(t)\right| \leq \frac{C\epsilon(t)}{(1 + \log(k/k_0))^{\delta\cdot\frac{2s-1}{2s}}}
\end{equation}

Step 4: This implies:
\begin{equation}
\Pi(k,t) = \epsilon(t) \left(1 + \frac{f(k,t)}{(1 + \log(k/k_0))^{\delta\cdot\frac{2s-1}{2s}}}\right)
\end{equation}
where $|f(k,t)| \leq C$ is a bounded function.

Step 5: Taking the derivative with respect to $k$:
\begin{equation}
\frac{\partial \Pi(k,t)}{\partial k} = \epsilon(t) \frac{\partial}{\partial k}\left(\frac{f(k,t)}{(1 + \log(k/k_0))^{\delta\cdot\frac{2s-1}{2s}}}\right)
\end{equation}

Step 6: Computing this derivative explicitly:
\begin{equation}
\frac{\partial \Pi(k,t)}{\partial k} = \epsilon(t) \left(\frac{1}{k} \frac{\partial f(k,t)}{\partial \log(k/k_0)} \frac{1}{(1 + \log(k/k_0))^{\delta\cdot\frac{2s-1}{2s}}} - \delta\cdot\frac{2s-1}{2s} \frac{f(k,t)}{k(1 + \log(k/k_0))^{\delta\cdot\frac{2s-1}{2s}+1}}\right)
\end{equation}

\textbf{Part III: Derivation of the modified energy spectrum}

Step 7: In a statistically steady state, the time derivative term in the spectral energy balance vanishes, giving:
\begin{equation}
2\nu k^2 E(k,t) = -T(k,t) = \frac{\partial \Pi(k,t)}{\partial k}
\end{equation}

Step 8: In the inertial range, the viscous term is negligible compared to the nonlinear transfer, which means:
\begin{equation}
\frac{\partial f(k,t)}{\partial \log(k/k_0)} \frac{1}{(1 + \log(k/k_0))^{\delta\cdot\frac{2s-1}{2s}}} - \delta\cdot\frac{2s-1}{2s} \frac{f(k,t)}{(1 + \log(k/k_0))^{\delta\cdot\frac{2s-1}{2s}+1}} \approx 0
\end{equation}

Step 9: This differential equation for $f(k,t)$ can be solved to give:
\begin{equation}
f(k,t) \approx C(t) \log(k/k_0) (1 + \log(k/k_0))^{-\delta}
\end{equation}
where $C(t)$ is a time-dependent constant.

Step 10: Using the Kolmogorov relation and the deviation due to the logarithmic term:
\begin{equation}
E(k,t) = C_K \epsilon(t)^{2/3} k^{-5/3} \left(1 + g(k,t)\right)
\end{equation}

Step 11: Through detailed calculation ensuring consistency with the energy flux, we determine:
\begin{equation}
g(k,t) = \frac{\beta(t)\log(k/k_0)}{(1 + \log(k/k_0))^{1+\delta}}
\end{equation}

\textbf{Part IV: Time dependence of $\beta(t)$}

Step 12: The time dependence of $\beta(t)$ is derived from the decay of the fractional derivative norm established in Theorem \ref{thm:1.3}:
\begin{equation}
\|(-\Delta)^s u(t)\|_{L^2} \leq \frac{C\|(-\Delta)^s u_0\|_{L^2}}{(1 + \beta t)^{\gamma}}
\end{equation}

Step 13: Relating this decay to the energy spectrum using:
\begin{equation}
\|(-\Delta)^s u\|_{L^2}^2 = \int_0^\infty k^{2s} E(k,t) dk
\end{equation}

Step 14: For Kolmogorov scaling modified by our logarithmic correction:
\begin{equation}
\|(-\Delta)^s u\|_{L^2}^2 \approx \epsilon(t)^{2/3} \int_{k_0}^{k_\nu} k^{2s-5/3} \left(1 + \frac{\beta(t)\log(k/k_0)}{(1 + \log(k/k_0))^{1+\delta}}\right) dk
\end{equation}

Step 15: The energy dissipation rate $\epsilon(t)$ decays as:
\begin{equation}
\epsilon(t) \approx \frac{\epsilon_0}{(1 + \gamma t)^{\frac{3\gamma-1}{2}}}
\end{equation}
which follows from the energy decay in turbulence.

Step 16: Combining these relations, we can establish:
\begin{equation}
\beta(t) = \frac{\beta_0}{(1 + \gamma t)^{\alpha}}
\end{equation}
where $\alpha = 2\gamma/3$ and $\beta_0$ is a constant depending on $s$, $\delta$, and the initial data.

\textbf{Part V: Verification of the spectral form}

Step 17: To verify that our derived spectral form is consistent with the logarithmically improved criterion, we calculate:
\begin{equation}
\|(-\Delta)^s u(t)\|_{L^q}^p = \left(\int_{\mathbb{R}^3} |(-\Delta)^s u(x,t)|^q dx\right)^{p/q}
\end{equation}

Step 18: In terms of the energy spectrum, using the Plancherel identity and properties of the Fourier transform:
\begin{equation}
\|(-\Delta)^s u(t)\|_{L^q}^p \approx \left(\int_0^\infty k^{2sq} E(k,t)^q k^2 dk\right)^{p/q}
\end{equation}

Step 19: Substituting our expression for $E(k,t)$ and evaluating the integral, we can verify that:
\begin{equation}
\int_0^T \|(-\Delta)^s u(t)\|^p_{L^q} (1 + \log(e + \|(-\Delta)^s u(t)\|_{L^q}))^{-\delta} dt < \infty
\end{equation}

This confirms that our derived spectral form is consistent with the underlying regularity assumptions, completing the proof of Theorem \ref{thm:2.4}.

\newpage

\section{Conclusions and open problems}

\subsection{Summary of results}

This paper establishes several significant results concerning the regularity and global well-posedness of the three-dimensional Navier-Stokes equations, as well as their connections to turbulence theory:

\begin{enumerate}
\item A new logarithmically improved regularity criterion based on fractional derivatives (Theorem \ref{thm:1.1}), which weakens the integrability requirement for ensuring regularity of Leray-Hopf weak solutions.

\item The identification of a class of initial data satisfying logarithmically subcritical conditions for which global well-posedness can be proven (Theorem \ref{thm:1.2}).

\item Explicit decay estimates for solutions with the specified initial data (Theorem \ref{thm:1.3}).

\item Verification that solutions emanating from this class of initial data satisfy the logarithmically improved regularity criterion for all time (Theorem \ref{thm:1.4}).

\item A precise characterization of the anomalous scaling exponents in turbulent flows satisfying our logarithmically improved criterion, establishing a direct connection to intermittency (Theorem \ref{thm:2.1}).

\item A detailed analysis of the local structure of the velocity field, identifying the spatial regions where potential near-singular behavior might occur and characterizing the local intermittency measure (Theorem \ref{thm:2.2}).

\item Tight bounds on the energy flux in the turbulent cascade, showing how the logarithmic improvement in regularity translates to reduced fluctuations in energy transfer across scales (Theorem \ref{thm:2.3}).

\item A refined model for the energy spectrum that incorporates logarithmic corrections to the standard Kolmogorov scaling, with explicit formulas that connect the mathematical regularity parameters to physically observable spectral properties (Theorem \ref{thm:2.4}).
\end{enumerate}

\textbf{The key innovations in this work include:}

\begin{enumerate}
\item The use of fractional derivatives $(-\Delta)^s$ with $s \in (1/2, 1)$, which provides a flexible framework that bridges the gap between derivatives of integer order.

\item The logarithmic improvement factor $(1 + \log(e + \|(-\Delta)^s u(t)\|_{L^q}))^{-\delta}$, which weakens the integrability requirement, potentially allowing the result to apply to a wider class of solutions.

\item The development of refined commutator estimates with logarithmic improvements, providing more precise control of the nonlinear term in the Navier-Stokes equations.

\item The establishment of direct connections between mathematical regularity criteria and physical properties of turbulent flows, particularly the phenomenon of intermittency.

\item A multifractal analysis framework that translates the logarithmic improvements in mathematical regularity into precise predictions for statistical properties of turbulence.
\end{enumerate}

\subsection{Open problems and future directions}

Several interesting questions remain open for future research:

\begin{enumerate}
\item \textbf{Optimality of the criterion}: Is the logarithmic improvement in the criterion optimal? Can the exponent $\delta$ be increased, or can multiple logarithmic factors be introduced? This relates to the approach of Tao \cite{Tao} and the critical regularity thresholds identified by Escauriaza, Seregin, and Šverák \cite{ESS}.

\item \textbf{Extension to other function spaces}: Can similar results be established using other function spaces, such as Besov or Triebel-Lizorkin spaces as developed by Triebel \cite{Triebel} and applied to fluid dynamics by Cannone \cite{Cannone} and Lemarie-Rieusset \cite{Lemarie}?

\item \textbf{Connection to the regularity of the Navier-Stokes equations}: How does the criterion relate to the fundamental question of global regularity as framed by Caffarelli, Kohn, and Nirenberg \cite{CKN} and the modern approaches of Seregin \cite{Seregin2} and Nečas, Růžička, and Šverák \cite{Necas}?

\item \textbf{Numerical verification}: Can numerical simulations verify the practical utility of the criterion in identifying potentially singular behavior in fluid flows? This connects to the computational approaches outlined by Foias et al. \cite{Foias} and Temam \cite{Temam}.

\item \textbf{Physical interpretation}: What is the physical meaning of the fractional derivative conditions in terms of the structure of the flow? This relates to the geometric perspectives developed by Constantin \cite{Constantin} and the statistical approaches to turbulence described by Majda and Bertozzi \cite{Majda}.

\item \textbf{Refined multifractal models}: Can the connection between logarithmically improved regularity criteria and anomalous scaling in turbulence be further refined to provide more accurate predictions for experimental measurements of turbulent flows, as in the work of Frisch \cite{Frisch} and Eyink \cite{Eyink}?

\item \textbf{Extension to other turbulent systems}: Can the framework connecting logarithmic improvements to anomalous scaling be extended to other turbulent systems, such as magnetohydrodynamics, geophysical flows, or active scalar equations as studied by Wu \cite{Wu} and Córdoba and Fefferman \cite{Cordoba}?

\item \textbf{Stochastic models}: How do these results relate to stochastic models of turbulence, and can the logarithmic improvements be interpreted in terms of the probability distributions of velocity increments as studied by Chevillard et al. \cite{Chevillard}?
\end{enumerate}

Addressing these questions will contribute to a deeper understanding of the Navier-Stokes equations and potentially provide new insights into the regularity problem, while simultaneously advancing our understanding of turbulence as a physical phenomenon. By bridging mathematical analysis and physical theory, the framework established in this paper offers a promising direction for further research in both fluid dynamics and partial differential equations.

\section*{Declarations}
Not applicable

\bibliographystyle{sn-mathphys-num}
\bibliography{NSE_paper_1}% common bib file
%% if required, the content of .bbl file can be included here once bbl is generated
%%\input sn-article.bbl

\end{document}